\documentclass[12pt]{amsart}
\usepackage{amssymb}
\usepackage{amsmath}
\usepackage{amsthm}
\usepackage{amstext}
\usepackage{amsopn}
\usepackage{mathrsfs}
\usepackage{latexsym}
\allowdisplaybreaks
\newtheorem{Definition}{Definition}[section]
\newtheorem{Proposition}{Proposition}[section]
\newtheorem{Lemma}{Lemma}[section]
\newtheorem{Theorem}{Theorem}[section]
\newtheorem{Corollary}{Corollary}[section]
\newtheorem{Remark}{Remark}[section]
\newtheorem{Example}{Example}[section]

\begin{document}
\bibliographystyle{plain}
\footnotetext{
This work is partially supported by MNiSW research grant No N N201 364436}
\title[Asymptotic properties of random matrices and pseudomatrices]
{Asymptotic properties of random \\matrices and pseudomatrices}
\author[R. Lenczewski]{Romuald Lenczewski}
\address{Romuald Lenczewski, \newline
Instytut Matematyki i Informatyki, Politechnika Wroc\l{}awska, \newline
Wybrze\.{z}e Wyspia\'{n}skiego 27, 50-370 Wroc{\l}aw, Poland  \vspace{10pt}}
\email{Romuald.Lenczewski@pwr.wroc.pl}
\maketitle

\begin{abstract}
We study the asymptotics of sums of matricially free random variables, 
called random pseudomatrices, and we compare it with that of random matrices
with block-identical variances. For objects of both types 
we find the limit joint distributions of blocks and give their Hilbert space 
realizations, using operators called `matricially free Gaussian operators'. 
In particular, if the variance matrices 
are symmetric, the asymptotics of symmetric blocks of random pseudomatrices agrees with 
that of symmetric random blocks.
We also show that blocks of random pseudomatrices are `asymptotically matricially free'
whereas the corresponding symmetric random blocks are `asymptotically symmetrically matricially free', where
symmetric matricial freeness is obtained from matricial freeness by an operation
of symmetrization. Finally, we show that row blocks of square, lower-block-triangular and block-diagonal
pseudomatrices are asymptotically free, monotone independent and
boolean independent, respectively. 
\end{abstract}

\section{Introduction and main results}
We have recently shown that the Hilbert space construction of
the free product of states on $C^{*}$-algebras given by Voiculescu [8]
can be generalized to a framework which exhibits some matricial features [5].

We considered arrays of noncommutative probability spaces, for which 
we defined the {\it matricially free product of states}.
The definition of this product is based on replacing the family of canonical 
unital *-representations of $C^{*}$-algebras on a free product of Hilbert spaces by a
diagonal-containing array of non-unital 
*-representations of $C^{*}$-algebras on 
the matricially free product of Hilbert spaces. 
The crucial point is that products of these 
representations imitate products of matrices rather than free products, although
the main features of the latter are still present.
We studied the associated concepts of noncommutative independence 
called {\it matricial freeness} and related to it {\it strong matricial freeness} which can be viewed as a scalar-type generalization of freeness underlying other fundamental notions of noncommutative independence (monotone and boolean)
and some of their generalizations (conditional freeness and conditional monotone independence).

More importantly, matricial freeness is closely related to random matrices. 
Their significance for free probability was discovered by Voiculescu [10], who showed that Gaussian random matrices with 
mutually independent entries are asymptotically free. This connection was later
generalized by Dykema [2] to non-Gaussian random matrices. 
Of course, the first indication that free probability might be related to random matrices was the central limit 
theorem for free random variables since its limit law was
the semicirle law obtained in the classical work of Wigner [13] as the limit distribution 
of certain self-adjoint random matrices with independent entries.

The main objects of interest in our study of matricial freeness and strong matricial freeness are diagonal-containing arrays of (in general, non-unital) subalgebras of a 
given unital algebra ${\mathcal A}$, 
$$
({\mathcal A}_{i,j})_{(i,j)\in J} 
\;\;\;{\rm with}\;\;\;\Delta\subseteq J\subseteq I\times I\;\;\;{\rm and}\;\;\;
\Delta=\{(j,j):j\in I\},
$$ 
equipped with an array of states $(\varphi_{i,j})$ on ${\mathcal A}$,
with respect to which notions of noncommutative independence are defined. 
Here, the states $\varphi_{i,j}$ and their kernels 
play a very similar role to that of one distinguished state $\varphi$ and its kernel 
in free probability. Apart from square arrays, of particular interest are lower- (upper-) triangular arrays related to monotone (anti-monotone) independence introduced by Muraki [6]. In all arrays, the roles of diagonal and off-diagonal entries are quite different, which is a characteristic feature of this theory reminding the random matrix theory. 

Let $(X_{i,j}(n))$ be an $n$-dimensional square array 
(or its diagonal-containing subarray) of self-adjoint random 
variables in a unital *-algebra ${\mathcal A}(n)$ 
which is matricially free with respect to the array 
$(\varphi_{i,j}(n))$ defined in terms 
of a family $(\varphi_j(n))_{1\leq j \leq n}$. The sums
$$
S(n)=\sum_{(i,j)\in J(n)}X_{i,j}(n),
$$
where $(J(n))$ is an appropriate sequence of sets,
called {\it random pseudomatrices}, remind random matrices, whereas
convex linear combinations
$$
\psi(n)=\frac{1}{n}\sum_{j=1}^{n}\varphi_{j}(n),
$$ 
where $n\in {\mathbb N}$, replace the states
$\tau(n)=\mathbb{E} \otimes {\rm tr(n)}$, classical expectation tensored with
normalized traces, under which distributions of random matrices are computed.

Under suitable assumptions on the distributions of the variables $X_{i,j}(n)$
in the states $\varphi_{i,j}(n)$, respectively, we can now study the asymptotic 
distributions of sums corresponding to blocks of random pseudomatrices in various states. 
These blocks are defined as sums of the form
$$
S_{p,q}(n)=\sum_{(i,j)\in N_{p}\times N_{q}}X_{i,j}(n)
$$
where $[n]:=\{1,2, \ldots , n\}=N_{1}\cup N_{2}\cup \ldots \cup N_{r}$ is a partition 
into disjoint subsets whose sizes grow proportionately to $n$ (it is convenient to think of intervals).
The dependence of the $N_{j}$'s on $n$ is supressed in the notation.
In particular, we assume that the variances $v_{i,j}(n):=\varphi_{i,j}(n)(X_{i,j}^{2}(n))$ 
are identical within blocks.

This leads to various asymptotic properties of sums of matricially 
free and strongly matricially free random variables. Let $\varphi_{i,j}(n)$
be defined according to
$$
\varphi_{j,j}(n)=\varphi(n)\;\;\;{\rm and}\;\;\; \varphi_{i,j}(n)=\varphi_{j}(n)\;\;\;
{\rm for}\;\;\;i\neq j,
$$
where $\varphi(n)$ is a distinguished state on ${\mathcal A}(n)$ 
and $(\varphi_{j}(n))_{1\leq j\leq n}$ is a family of additional states 
called {\it conditions} associated with $\varphi(n)$ 
for any $n\in {\mathbb N}$. Then, as $n\rightarrow \infty$, we obtain 
the following properties:
\begin{enumerate}
\item[(i)]
if the variances of square arrays are identical within blocks and rows, then the sums
$$
A_{p}(n):=\sum_{q=1}^{n}S_{p,q}(n), \;\;\;{\rm where}\;\;1 \leq p\leq r,
$$ 
are asymptotically free with respect to $\varphi(n)$,
\item[(ii)]
if the variances of block-lower-triangular arrays are identical within blocks and rows, 
then the sums 
$$
B_{p}(n):=\sum_{q=1}^{p}S_{p,q}(n),\;\;\;{\rm where}\;\; 1 \leq p\leq r,
$$ 
are asymptotically monotone independent with respect to $\varphi(n)$,
\item[(iii)]
if the variances of block-diagonal arrays are identical within blocks, then the sums
$C_{p}(n):=S_{p,p}(n)$, where $1\leq p\leq r$, are asymptotically boolean independent with respect to $\varphi(n)$.
\end{enumerate}
Let us add that non-asymptotic analogs of (i)-(ii) for finite sums of {\it strongly} 
matricially free random variables were proved in [5] and that (iii) also holds for finite sums 
of both matricially free and strongly matricially free random variables.
Recall that the `strongly matricially free product of states', 
on which the definition of strong matricial freeness is based, is obtained from the 
matricially free product by restriction. However, it is worth to remark that 
the difference between blocks of matricially free and strongly matricially free 
random variables disappears asymptotically.

The scheme of matricial freeness and its symmetrized version, 
in which ordered pairs are replaced by (non-ordered) sets consisting of one or two elements, called {\it symmetric matricial freeness}, is also used for arrays 
$(\varphi_{i,j}(n))$ defined in terms of a family $(\varphi_j(n))_{1\leq j\leq n}$ 
of states on ${\mathcal A}(n)$ according to
$$
\varphi_{i,j}(n)=\varphi_j(n)\;\;\;{\rm for}\;{\rm any}\;i,j,
$$
and any $n\in {\mathbb N}$. Here, we do not a priori assume 
that the states $\varphi_j(n)$ are conditions associated 
with a distinguished state. In the study of asymptotics, 
of importance become arrays defined by `normalized partial traces'
$$
\psi_{q}(n)=\frac{1}{n_{q}}\sum_{j\in N_{q}}\varphi_{j}(n)
$$ 
according to $\psi_{p,q}(n)=\psi_{q}(n)$
for any $p,q\in [r]$ and $n\in {\mathbb N}$, 
where $n_{q}$ is the cardinality of $N_{q}$. 
Then, as $n\rightarrow \infty$, we obtain the following properties:
\begin{enumerate}
\item[(iv)]
if the variances are identical within blocks, then the array $(S_{p,q}(n))$ is
asymptotically matricially free with respect to $(\psi_{p,q}(n))$,
\item[(v)]
if the variances are identical within blocks and form symmetric matrices, 
then the array $(Z_{p,q}(n))$ of symmetric blocks given by
$$
Z_{p,q}(n):=
\sum_{(i,j)\in N_{p,q}}X_{i,j}(n)
$$ 
where $N_{p,q}=(N_{p}\times N_{q})\cup (N_{q}\times N_{p})$,
is asymptotically symmetrically matricially free
with respect to $(\psi_{p,q}(n))$,
\item[(vi)]
if the variances are identical within blocks of symmetric random matrices,
then the array $(T_{p,q}(n))$ of symmetric random blocks
is asymptotically symmetrically matricially free with respect to $(\tau_{p,q}(n))$,
\end{enumerate}
where the array $(\tau_{p,q}(n))$ is defined by the family 
of partial traces $\tau_q=\mathbb{E} \otimes {\rm tr}_{q}(n)$
and ${\rm tr}_{q}(n)$ denotes the normalized trace over the 
vectors indexed by $N_{q}$. Here, by symmetric random blocks we understand symmetric blocks of symmetric random matrices in the approach of Voiculescu [10] and Dykema [2].
 
Concerning a non-asymptotic analog of (iv), 
we show in this paper that blocks $(S_{p,q}(n))$ are matricially free 
for any finite $n$ under the stronger assumption that the variables have 
block-identical distributions.
In turn, properties (v) and (vi) do not seem to have their non-asymptotic analogs.
Nevertheless, they are the reason why we call the sums $S(n)$ random pseudomatrices. 
On the other hand, let us also remark that we do not have an analog of (iv) 
for random matrices. Therefore, informally, random pseudomatrices can be viewed as objects which remind random matrices for large $n$ if we consider blocks with symmetric variances, but exhibit different features if the variances are not symmetric, which leads to triangular arrays and relations to monotone independence.   

In that connection we would like to mention the result of Shlyakhtenko [7] 
who established a connection between a class of symmetric random matrices called 
random band matrices and an operator-type generalization of freeness called
freeness with amalgamation. Our result (vi) gives a connection between another class of symmetric random matrices called symmetric random blocks with block-identical distributions and a scalar-type generalization of freeness called 
symmetric matricial freeness. At present, it is not clear to us whether there is a relation
between these two approaches. 

In our previous work, we expressed the limit distributions of random pseudomatrices with respect to
$\varphi(n)$'s and $\psi(n)$'s in terms of certain functions on 
the class of colored non-crossing pair partitions, as well as in terms of some 
`continued multifractions' ]5]. Both these realizations showed that the limit laws can be viewed as matricial
generalizations of the semicircle laws. Moreover, these distributions turned out to 
be related to s-free additive convolutions and the associated notion of 
freeness with subordination, or s-freeness [3], concepts motivated by the results of Voiculescu [11] and Biane [1] 
on analytic subordination in free probability. Using s-freeness, we also studied s-free and free multiplicative convolutions, which 
allowed us to establish relations between these convolutions (as well as some other multiplicative convolutions)
and certain classes of walks on appropriately defined products of graphs [4].

In this paper, we give Hilbert-space realizations of the limit distributions 
of random pseudomatrices and their blocks under $\varphi(n)$, $\psi_k(n)$
and $\psi(n)$, where the underlying Hilbert space in all considered cases is
the matricially free product 
$$
({\mathcal F}, \xi)=*_{i,j}^{M}({\mathcal F}_{p,q})
$$ 
of the $r$-dimensional array $({\mathcal F}_{p,q})$ of Fock spaces, where
the diagonal and off-diagonal Fock spaces are, respectively, 
free and boolean Fock spaces over one-dimensional Hilbert spaces. 
Using realizations on ${\mathcal F}$ and operators which play the role of
Gaussian operators in our approach, as well as their `symmetrizations',
we prove (i)-(vi). 

In Section 2, we present the combinatorics of colored non-crossing pair-partitions. 
In Section 3 we recall the basic notions and facts concerning 
matricially free random variables.
In Section 4 we introduce matricially free Gaussian operators on ${\mathcal F}$.
Using them, we find $\mathcal{F}$-realizations of the limit distributions 
for random pseudomatrices in Section 5.
In Section 6 we prove that blocks of (strongly) matricially free arrays of random 
variables with block-identical distributions are matricially free with respect to a 
suitably defined array of states. 
In turn, asymptotic matricial freeness of blocks of random pseudomatrices is proved 
in Section 7, where we also find ${\mathcal F}$-realizations of their limit joint distributions.
In Section 8 we introduce the notion of `symmetric matricial freeness'. 
In Section 9 we find the limit joint distributions of symmetric random blocks 
and their ${\mathcal F}$-realizations and we show that
symmetric random blocks are asymptotically symmetrically matricially free. 
In Section 10 we obtain results on asymptotic freeness and asymptotic monotone independence 
of rows of pseudomatrices.

\section{Combinatorics}

The combinatorics of our model is based on the class of {\it colored non-crossing pair partitions}, 
to which we assign certain products of matrix elements whose indices depend on the colorings of their blocks. 

For a given non-crossing pair partition $\pi$,
we denote by $\mathcal{B}(\pi)$, $\mathcal{L}(\pi)$ and $\mathcal{R}(\pi)$ the sets of its blocks, 
their left and right legs, respectively. If $\pi_{i}=\{l(i),r(i)\}$ and $\pi_{j}=\{l(j),r(j)\}$
are blocks of $\pi$ with left legs $l(i)$ and $l(j)$ and right legs $r(i)$ and $r(j)$, respectively, then
$\pi_i$ is {\it inner} with respect to $\pi_j$ if $l(j)<l(i)<r(i)<r(j)$.
In that case $\pi_j$ is {\it outer} with respect to $\pi_i$. 
It is the {\it nearest} outer block of $\pi_i$ if there is no block $\pi_k=\{l(k),r(k)\}$
such that $l(j)<l(k)<l(i)<r(i)<r(k)<r(j)$. Since the nearest outer block, if it exists, is unique,
we can write in this case $\pi_j=o(\pi_i)$, $l(j)=o(l(i))$ and $r(j)=o(r(i))$. 
If $\pi_i$ does not have an outer block, it is called a {\it covering} block.
It is convenient to extend each partition $\pi\in \mathcal{NC}_{m}^{2}$
to the partition $\widehat{\pi}$ obtained from $\pi$ by adding one block, say $\pi_0=\{0, m+1\}$,
called the {\it imaginary block}. 

Let $F_{r}(\pi)$ be the set of all mappings 
$f:{\mathcal B}(\pi)\rightarrow [r]$ called {\it colorings} 
of the blocks of $\pi$ by the set $[r]:=\{1,2, \ldots ,r\}$. 
Then the pair $(\pi,f)$ plays the role
of a colored partition. Its blocks will be denoted 
$$
{\mathcal B}(\pi,f)=\{(\pi_1,f), 
(\pi_2,f), 
\ldots , (\pi_k,f)\},
$$ 
where we use pairs $(\pi_i,f)$ since to each block $\pi_i$ we shall assign entries of a 
matrix which depend on the colors of both $\pi_i$ and $o(\pi_i)$ for any $i\in [k]$.
If the imaginary block is used, it is convenient to assume that 
it is also colored by a number from the set $[r]$. 
In Fig. 1 we give an example with the notions defined above.

\begin{Definition}
{\rm Let $(\pi,f)$ be a colored non-crossing partition with blocks
as above, where $f\in F_{r}(\pi)$ and let $B\in M_{r}({\mathbb R})$ 
be given. For any $0\leq j \leq r$ we define 
$$
b_{j}(\pi, f)=b_{j}(\pi_{1},f)
b_{j}(\pi_2,f)
\ldots b_{j}(\pi_k,f)
$$
where the functions $b_j:{\mathcal B}(\pi,f)\rightarrow {\mathbb R}$ are
given by the following rules:
\begin{enumerate} 
\item
$b_{j}(\pi_i,f)=b_{p,q}$ if $f(\pi_i)=p$ and $f(o(\pi_i))=q$, where
$0\leq j\leq r$,
\item  
$b_j(\pi_i,f)=b_{p,j}$ if $f(\pi_i)=p$ and $\pi_i$ does not 
have outer blocks, where $1\leq j \leq r$,
\item 
$b_{0}(\pi_i,f)=b_{p,p}$ if $f(\pi_i)=p$ and 
$\pi_i$ does not have outer blocks.
\end{enumerate}}
\end{Definition}

The index $j\in [r]$ in $b_{j}(\pi,f)$ can be 
interpreted as the color of the imaginary block (if $j=0$, then 
the imaginary block is not needed).
Finally, $\mathcal{NC}_{m}^{2}=\emptyset$ for $m$ odd and thus we shall 
understand in this case that the 
summation over $\pi\in \mathcal{NC}_{m}^{2}$ gives zero.

When we sum these products over all possible colorings, we obtain numbers
$$
b_{j}(\pi)=\sum_{f\in F_{r}(\pi)}b_{j}(\pi,f)
$$
for any $0\leq j \leq r$, which are used to express the limit 
laws of random pseudomatrices. 

\begin{Example}
{\rm If we consider the partition $\pi$ given by the diagram in Fig.1, we obtain
$$
b_{0}(\pi)=\sum_{i,k,l}b_{i,l}b_{k,l}b_{l,l}
\;\;\;{\rm and}\;\;\;
b_{j}(\pi)=\sum_{i,k,l}b_{i,l}b_{k,l}b_{l,j}
$$
where $j\in [r]$ and all indices run over $[r]$, which corresponds to all possible 
colorings from $F_{r}(\pi)$ and is adequate for a square array $B\in M_{r}({\mathbb R})$.
For other arrays, we obtain the same formulas except that range of the summation is smaller.}
\end{Example}

Entries of the matrix $B$ will be related to variances of matricially free random variables
obtained when computing their mixed moments.
Therefore, it is natural to begin with an array
$(a_{i,j})$ of random variables and assign 
a variable to each leg of the considered partition.
If the array is square, we 
assume that $i,j\in [r]$, whereas in other cases the pairs $(i,j)$ belong
to a proper subset $J\subset [r]\times [r]$ which includes the diagonal.
To a given non-crossing partition $\pi\in \mathcal{NC}_{m}^2$, where
$m$ is even, we can now associate a product of these variables in which 
indices can be interpreted as colors taken from the set $[r]$. These indices
are related to each other in a natural way which refers to
the way matricially free random variables can be multiplied to
give a non-trivial contribution to the limit laws. The definition given below
is based on this relation. 

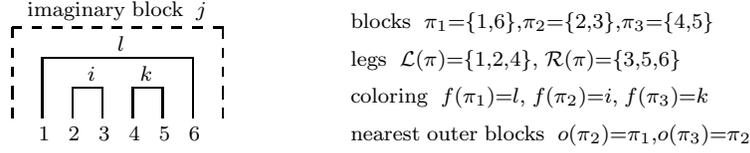
\begin{figure}
\unitlength=1mm
\special{em:linewidth 0.4pt}
\linethickness{0.4pt}
\begin{picture}(120.00,30.00)(52.00,5.00)
\put(69.00,10.00){\line(0,1){8.00}}
\put(73.00,10.00){\line(0,1){4.00}}
\put(77.00,10.00){\line(0,1){4.00}}
\put(81.00,10.00){\line(0,1){4.00}}
\put(85.00,10.00){\line(0,1){4.00}}
\put(89.00,10.00){\line(0,1){8.00}}

\put(65.00,10.00){\line(0,1){1.25}}
\put(65.00,12.50){\line(0,1){1.25}}
\put(65.00,13.00){\line(0,1){1.25}}
\put(65.00,15.50){\line(0,1){1.25}}
\put(65.00,18.00){\line(0,1){1.25}}
\put(65.00,20.75){\line(0,1){1.25}}

\put(93.00,10.00){\line(0,1){1.25}}
\put(93.00,12.50){\line(0,1){1.25}}
\put(93.00,13.00){\line(0,1){1.25}}
\put(93.00,15.50){\line(0,1){1.25}}
\put(93.00,18.00){\line(0,1){1.25}}
\put(93.00,20.60){\line(0,1){1.25}}

\put(68.50,7.00){$\scriptstyle{1}$}
\put(72.50,7.00){$\scriptstyle{2}$}
\put(76.50,7.00){$\scriptstyle{3}$}
\put(80.50,7.00){$\scriptstyle{4}$}
\put(84.50,7.00){$\scriptstyle{5}$}
\put(88.50,7.00){$\scriptstyle{6}$}

\put(75.00,15.00){$\scriptstyle{i}$}
\put(82.00,15.00){$\scriptstyle{k}$}
\put(79.00,19.00){$\scriptstyle{l}$}
\put(67.00,23.50){$\scriptstyle{{\rm imaginary\;block}\;\;j}$}

\put(65.00,22.00){\line(1,0){1.25}}
\put(67.67,22.00){\line(1,0){1.25}}
\put(70.33,22.00){\line(1,0){1.25}}
\put(72.98,22.00){\line(1,0){1.25}}
\put(75.63,22.00){\line(1,0){1.25}}
\put(78.30,22.00){\line(1,0){1.25}}
\put(80.95,22.00){\line(1,0){1.25}}
\put(83.60,22.00){\line(1,0){1.25}}
\put(86.25,22.00){\line(1,0){1.25}}
\put(89.00,22.00){\line(1,0){1.25}}
\put(91.70,22.00){\line(1,0){1.30}}

\put(69.00,18.00){\line(1,0){20.00}}
\put(73.00,14.00){\line(1,0){4.00}}
\put(81.00,14.00){\line(1,0){4.00}}
\put(110.00,22.00){${\scriptstyle{\rm blocks}\;\;\pi_1=\{1,6\}, \pi_2=\{2,3\}, \pi_3=\{4,5\}}$}
\put(110.00,17.00){${\scriptstyle{\rm legs}\;\;{\mathcal L}(\pi)=\{1,2,4\},\;{\mathcal R}(\pi)=\{3,5,6\}}$}
\put(110.00,12.00){${\scriptstyle{\rm coloring}\;\;f(\pi_1)=l,\; f(\pi_2)=i,\; f(\pi_3)=k}$}
\put(110.00,7.00){${\scriptstyle{\rm nearest\;outer\;blocks}\;\;o(\pi_2)=\pi_1, o(\pi_3)=\pi_2}$}

\end{picture}
\caption{A colored non-crossing partition}
\end{figure}
\begin{Definition}
{\rm We will say that the partition $\pi \in \mathcal{NC}_{m}^{2}$ is {\it adapted}
to a tuple of numbers $(p_1,q_1, \ldots, p_m,q_m)$ from the set $[r]$ if
\begin{enumerate}
\item
$(p_i,q_i)=(p_j,q_j)$ whenever $\{i,j\}$ is a block of $\pi$,
\item
$q_k=p_{o(k)}$ whenever $k\in \mathcal{R}(\pi)$.
\end{enumerate}
The set of such partitions will be denoted $\mathcal{NC}_{m}^{2}(p_1,q_1, \ldots , p_m,q_m)$.}
\end{Definition}
If $\pi\in\mathcal{NC}_{m}^{2}(p_1,q_1,\ldots, p_m,q_m)$, then 
the tuple $(p_1,q_1, \ldots , p_m,q_m)$ defines a unique coloring of 
both $\pi$ and $\widehat{\pi}$ 
in which the block containing $k\in \mathcal{L}(\pi)$ is colored by $p_k$ for any $k$ 
and the imaginary block is colored by $q_m$.

\begin{Example}
{\rm Consider the partition $\pi$ shown in Fig.1 and let $p_1,q_1, \ldots , p_6,q_6\in [6]$ be given.
Conditions (1)-(2) of Definition 2.2 lead to equations
$p_1=p_6,p_2=p_3,p_4=p_5$ and $q_2=q_3=p_1$, $q_4=q_5=p_1$, $q_1=q_6$, which 
gives four independent colors. Setting $p_1=l , p_2=i, p_4=k $ and $q_6=j$, 
we obtain the coloring of Fig.1.
Note in this context that $p_1,p_2,p_4$ may be interpreted as colors
of the left legs of $\pi$ and $q_6$ as the color of the imaginary block.
These equations lead to a product of variables of the form
$$
a_{p_{1},q_6}a_{p_2,p_1}a_{p_2,p_1}a_{p_4,p_1}a_{p_4,p_1}a_{p_1,q_6},
$$
where each variable is taken from an array $(a_{p,q})$ of
matricially free random variables and is associated with one inner-outer pair 
of blocks of $\widehat{\pi}$.  
}
\end{Example}

We close this section with the symmetrized version of Definition 2.2 which will 
be needed in a combinatorial formula for the limit joint distribution of
symmetric random blocks. The symmetrization is obtained by replacing ordered
pairs by subsets.

\begin{Definition}
{\rm We will say that the partition $\pi \in \mathcal{NC}_{m}^{2}$ is {\it adapted} to 
a tuple of subsets $(\{p_1,q_1\}, \ldots , \{p_m,q_m\})$ of the set $[r]$ if
\begin{enumerate}
\item
$\{p_i,q_i\}=\{p_j,q_j\}$ whenever $\{i,j\}$ is a block of $\pi$,
\item
$\{p_k,q_k\}\cap \{p_{o(k)},q_{o(k)}\}\neq \emptyset$ 
whenever $k\in \mathcal{R}(\pi)$.
\end{enumerate}
The set of such partitions will be denoted
$\mathcal{NC}_{m}^{2}(\{p_1,q_1\},\ldots,\{p_m,q_m\})$.}
\end{Definition}
If $\pi\in\mathcal{NC}_{m}^{2}(\{p_1,q_1\},\ldots, \{p_m,q_m\})$,
the tuple $(\{p_1,q_1\}, \ldots , \{p_m,q_m\})$ defines colorings of 
both $\pi$ and $\widehat{\pi}$, called {\it admissible},
in which each block containing $k\in \mathcal{R}(\pi)$ is colored by $p_k$ or $q_k$
and the imaginary block is colored by $p_m$ or $q_m$.

\begin{Example}
{\rm Consider again the partition shown in Fig.1 
and let $p_1,q_1, \ldots , p_6,q_6\in [6]$ be given.
Conditions (1)-(2) of Definition 2.3 take the form 
\begin{enumerate}
\item 
$\{p_1,q_1\}=\{p_6,q_6\}$, $\{p_2,q_2\}=\{p_3,q_3\}$
and $\{p_4,q_4\}=\{p_5,q_5\}$,
\item
$\{p_3,q_3\}\cap \{p_1,q_1\}\neq \emptyset$
and $\{p_5,q_5\}\cap \{p_1,q_1\}\neq \emptyset$.
\end{enumerate}
For instance, they are satisfied if 
$(p_1,q_1)=(q_6,p_6)$, $(p_2,q_2)=(q_3,p_3)$, $(p_4,q_4)=(q_5,p_5)$,
with $q_i=p_{i+1}$ for any $i=1, \ldots , 5$,
which corresponds to the product of variables of the form
$$
a_{p_1,p_6}a_{p_6,p_3}a_{p_3,p_6}a_{p_6,p_5}a_{p_5,p_6}a_{p_6,p_1}, 
$$
where we have four independent indices $p_1,p_3,p_5,p_6$, among which $p_1$ colors the imaginary block 
and $p_3,p_5,p_6$ color the blocks of $\pi$ (they are associated with the right legs of $\pi$).
If we replace some indices from the set $\{p_1,p_3,p_5,p_6\}$ by the corresponding $q_i$'s, we 
obtain other admissible colorings and the associated products of variables.}
\end{Example}

\section{Matricially free random variables} 

Let us recall the definition of matricially free random variables as well as
the main results of [5], where we refer the reader for details.

Let ${\mathcal A}$ be a unital algebra with an array $({\mathcal A}_{i,j})$ of not necessarily unital subalgebras of 
${\mathcal A}$ and let $(\varphi_{i,j})$ be a family of states on ${\mathcal A}$.  
Here, by a state on ${\mathcal A}$ we understand a normalized linear functional.
Further, we assume that each ${\mathcal A}_{i,j}$ has an {\it internal unit} $1_{i,j}$, for which
it holds that $a1_{i,j}=1_{i,j}a=a$ for any $a\in {\mathcal A}_{i,j}$, and
that the unital subalgebra ${\mathcal I}$ of ${\mathcal A}$ generated by all internal 
units is commutative. If ${\mathcal A}$ is a unital *-algebra, then, in addition, we require that 
the considered functionals are positive, the subalgebras are 
*-subalgebras and the internal units are projections. 

However, the definitions given below are slightly more general 
than those in [5]. Namely, in contrast to the formulation given there, 
we do not assume any particular form of the considered array $(\varphi_{i,j})$. 
Instead, we will distinguish the {\it diagonal} states as those of the form 
$\varphi_{j,j}$, where $j\in I$, but we will not assume that they all coincide. This will enable us to use the concept of 
matricial freeness for a wider class of arrays, which turns out convenient in the formulation 
of our results.

We shall use the subsets of $(I\times I)^{m}$ of the form
$$
\Lambda_{m}=\{((i_1,i_2), (i_2,i_3), \ldots , (i_m,i_{m+1})): (i_1,i_2)\neq (i_2,i_3)\neq \ldots \neq (i_m,i_{m+1})\},
$$ 
where $m\in {\mathbb N}$, with their union denoted $\Lambda=\bigcup_{m=1}^{\infty}\Lambda_m$. 
\begin{Definition}
{\rm We say that $(1_{i,j})$ is a {\it matricially free array of units} 
associated with $({\mathcal A}_{i,j})$ and $(\varphi_{i,j})$ if for any 
diagonal state $\varphi$ it holds that
\begin{enumerate}
\item
$\varphi(u_1au_2)=\varphi(u_1)\varphi(a)\varphi(u_2)$ 
for any $a\in {\mathcal A}$ and $u_1,u_2\in {\mathcal I}$,
\item
if $a_{k}\in {\mathcal A}_{i_k,j_k}\cap {\rm Ker}\varphi_{i_k,j_k}$, where $1<k\leq m$, then
$$
\varphi(a1_{i_1,j_1}a_{2}\ldots a_m)=
\left\{
\begin{array}{cc}
\varphi(aa_{2} \ldots a_n) & {\rm if}\;\;((i_{1},j_{1}), \ldots , (i_{m},j_m))\in \Lambda\\
0 & {\rm otherwise}
\end{array}
\right..
$$
where $a\in {\mathcal A}$ is arbitrary and $(i_1,j_1)\neq \ldots \neq (i_m,j_m)$.
\end{enumerate}}
\end{Definition}
\begin{Definition}
{\rm We say that $({\mathcal A}_{i,j})$ is 
{\it matricially free} with respect to $(\varphi_{i,j})$ if  
\begin{enumerate}
\item for any $a_{k}\in {\rm Ker}\varphi_{i_k,j_k}\cap {\mathcal A}_{i_k,j_k}$, 
where $k\in [m]$ and $(i_1,j_1)\neq \ldots \neq (i_m,j_m)$, and for any diagonal state 
$\varphi$, it holds that
$$
\varphi(a_1a_2\ldots a_m)=0
$$
\item
$(1_{i,j})$ is a matricially free array of units associated with 
$({\mathcal A}_{i,j})$ and $(\varphi_{i,j})$.
\end{enumerate}}
\end{Definition}

\begin{Definition}
{\rm  
The array of variables $(a_{i,j})$ in a unital algebra (*-algebra) ${\mathcal A}$ will be called 
{\it matricially free} (*-{\it matricially free}) with respect to $(\varphi_{i,j})$ if
there exists an array of elements (projections) $(1_{i,j})$ which is a 
matricially free array of units associated with ${\mathcal A}$ and $(\varphi_{i,j})$ 
and such that the array of algebras (*-algebras) generated by $a_{i,j}$ and $1_{i,j}$, 
respectively, is matricially free with respect to $(\varphi_{i,j})$.
}
\end{Definition}

Slightly less general was the setting given in [5], 
where we assumed that all diagonal states coincide with some 
distinguished state $\varphi$ and the off-diagonals states agree
with a family of additional states $(\varphi_j)_{j\in I}$,
called {\it conditions} associated with $\varphi$, 
which are defined by 
$$
\varphi_{j}(a)=\varphi(c_{j}ab_{j})
$$
for some $c_j,b_j\in {\mathcal A}_{j,j}\cap {\rm Ker}\varphi$ such that
$\varphi(c_jb_j)=1$ (if ${\mathcal A}$ is a *-algebra, $c_{j}=b_{j}^{*}$).
In this setting, it is natural to assume that $\varphi$ and $\varphi_{j}$'s 
are normalized according to
$$
\varphi(1_{i,j})=\delta_{i,j}\;\;\; {\rm and}\;\;\; \varphi_{j}(1_{i,k})=\delta_{j,k}
$$ 
for any $i,j,k$. However, in this paper we will also use other arrays, for instance such 
in which all states in the $j$-th column agree with some $\varphi_{j}$, 
where $(\varphi_{j})_{j\in I}$ is a given family of states on ${\mathcal A}$
which are not a priori assumed to be conditions associated with some distinguished state.
This motivates the following definition.

\begin{Definition}
{\rm Let $\varphi$ and $\varphi_{j}, j\in I$, be states on ${\mathcal A}$ and let 
$(\varphi_{i,j})$, where $(i,j)\in J$
and $\Delta\subseteq J \subseteq I\times I$, be an array of states on ${\mathcal A}$. 
If $\varphi_{j,j}=\varphi$ and $\varphi_{i,j}=\varphi_{j}$ for any $(i,j)\in J$, we will say 
$(\varphi_{i,j})$ is {\it defined by $\varphi$ and the family $(\varphi_{j})_{j\,\in I}$}. 
In turn, if $\varphi_{i,j}=\varphi_j$ for any $(i,j)\in J$, we will say that 
$(\varphi_{i,j})$ is {\it defined by the family $(\varphi_j)_{j\in I}$}.}
\end{Definition}

In this context, let us remark that a Hilbert space setting, similar 
to that in [5], can be given for a family $(\varphi_{j})_{j\in I}$ of product 
states on the unital *-algebra 
$$
{\mathcal A}:=\sqcup_{i,j}{\mathcal A}_{i,j}
$$ 
which are not defined as conditions associated with some distinguished state 
on ${\mathcal A}$. Here, $\sqcup_{i,j}{\mathcal A}_{i,j}$ stands for the free product 
of $C^*$-algebras without identification of units which is assumed to contain 
the empty word playing the role of the algebra unit. 

In this `tracial' framework, if $({\mathcal H}_{i,j},\pi_{i,j},\xi_{i,j})$ is the
array of GNS triples associated with an array 
$({\mathcal A}_{i,j},\varphi_{i,j})$ of noncommutative probability spaces,
the underlying product Hilbert spaces are of the form
$$
{\mathcal H}_{j}=\left({\mathbb C}\xi_{j,j}\oplus \bigoplus_{m=1}^{\infty}\bigoplus_{(i_1,i_2)\neq \ldots \neq (i_{m-1},j)\neq (j,j)}
{\mathcal H}_{i_1,i_2}^{0}\otimes \ldots \otimes {\mathcal H}_{j,j}^{0}\right)
$$
where $j\in I$ and ${\mathcal H}_{i,j}^{0}$ is the orthocomplement of ${\mathbb C}\xi_{i,j}$ in 
${\mathcal H}_{i,j}$ for each $i,j$, and their direct sum is 
$$
{\mathcal H}=\bigoplus_{j\in I}{\mathcal H}_{j},
$$ 
with the canonical inner product used in all direct sums.
The space ${\mathcal H}$ replaces the matricially free product 
of Hilbert spaces $*_{i,j}^{M}({\mathcal H}_{i,j},\xi_{i,j})$ used in [5].

Nevertheless, ${\mathcal H}$ can be embedded in 
a larger matricially free product of Hilbert spaces, namely 
$({\mathcal H}', \xi)=*_{i,j}^{M}({\mathcal H}_{i,j}', \xi_{i,j}')$, where
$$
{\mathcal H}_{i,j}'=\left\{
\begin{array}{cl}
{\mathcal H}_{i,j}& {\rm if} \;\;i\neq j\\
{\mathcal H}_{j,j}\oplus {\mathbb C}\xi_{j,j}'& {\rm if}\;\; i=j
\end{array}
\right.
\;\;\;\;{\rm and}\;\;\;\;
\xi_{i,j}'=\left\{
\begin{array}{ll}
\xi_{i,j} & {\rm if} \;\;i\neq j\\
\xi_{j,j}' & {\rm if}\;\; i=j
\end{array}
\right.
$$
and $\xi_{j,j}'$ is an additional unit vector for each $j\in I$.
In order to define appropriate product states on ${\mathcal A}$ 
associated with vectors $\xi_{j,j}$, we first trivially extend each cyclic representation 
$\pi_{j,j}: {\mathcal A}_{j,j}\rightarrow B({\mathcal H}_{j,j})$ 
to a non-cyclic representation $\pi_{j,j}':{\mathcal A}_{j,j}\rightarrow B({\mathcal H}_{j,j}')$,
and we keep the off-diagonal representations unchanged, namely $\pi_{i,j}'=\pi_{i,j}$ for $i\neq j$.  
Then, we take the free product $\lambda'=*_{i,j}^{M}\pi_{i,j}'$ and observe that 
${\mathcal H}$ is left invariant by $\lambda'(a)$ for any $a\in {\mathcal A}$. 
This allows us to define the unital *-representation 
$$
\lambda: {\mathcal A} \rightarrow B({\mathcal H})\;\;\;\;{\rm as}\;\;\;\;\lambda=\lambda'|{\mathcal H}
$$
which then leads to product states $\varphi_{j}$ on ${\mathcal A}$ 
defined by 
$$
\varphi_{j}= \eta_{j}\circ \lambda,
$$
where $\eta_{j}$ is the vector state on $B({\mathcal H})$ associated with  
$\xi_{j,j}$, where $j\in I$. 
An equivalent formulation in terms of partial isometries 
leading to the product states $\varphi_{j}$ can also be given.

The proposition given below provides the main motivation for 
allowing a more general class of arrays in the definition of 
matricial freeness.
\begin{Proposition}
The array of $C^{*}$-algebras $({\mathcal A}_{i,j})$ viewed as *-subalgebras of 
${\mathcal A}=\sqcup_{i,j}{\mathcal A}_{i,j}$ is matricially free with respect 
to the array $(\psi_{i,j})$ of states on ${\mathcal A}$ defined by 
the family $(\varphi_{j})_{j\in I}$ introduced above.    
\end{Proposition}
{\it Proof.}
In fact, it can be seen that an analog of [5, Proposition 2.3] holds for 
the states $\varphi_j:=\eta_{j}\circ \lambda$ defined above. However,
in this case it suffices to take $a_{k}\in {\mathcal A}_{i_k,j_k}$, 
where $k \in [n]$ and $(i_1,j_1)\neq \ldots \neq (i_n,j_n)$ (thus, 
we do not assume that $(i_1,j_1)\neq (j,j)\neq (i_n,j_n)$). 
If $\,a_{k}\in {\rm Ker}\,\varphi_{i_k,j_k}$ for $k \in [n]$, then
$$
\varphi_{j}(a_1a_2\ldots a_n)=0
$$ 
for each $j$. Moreover, if $\,a_r=1_{i_r,j_r}$ and $a_{m}\in {\rm Ker}\,\varphi_{i_m,j_m}$ for $r < m \leq n$, then
$$
\varphi_{j}(a_1\ldots a_n)=\left\{
\begin{array}{cc}
\varphi_{j}(a_1\ldots a_{r-1}a_{r+1}\ldots a_n)&\;{\rm if}\;\; ((i_r,j_r), \ldots , (i_n,j_n))\in \Lambda\\
0 & {\rm otherwise}
\end{array}\right..
$$
Finally, for any $a\in {\mathcal A}$, $u_1,u_2\in {\mathcal I}$ and $i,j,k\in I$, it holds that 
$$
\varphi_{j}(u_1au_2)=
\varphi_{j}(u_1)\varphi_{j}(a)\varphi_{j}(u_2)\;\;{\rm and}\;\;
\varphi_{j}(1_{i,k})=\delta_{j,k}.
$$
All these facts imply that $({\mathcal A}_{i,j})$ is matricially free with respect 
to the array $(\psi_{i,j})$, where $\psi_{i,j}=\varphi_{j}$ for any $i,j$. 
\hfill$\blacksquare$

\begin{Remark}
{\rm 
In this context, let us observe that if $({\mathcal A}_{i,j})$ is matricially free
with respect to the array defined by $\varphi$ and the associated conditions 
$(\varphi_{j})_{j\in I}$, then $({\mathcal A}_{i,j})$ is not, in general, matricially free with respect 
to the array defined by $(\varphi_j)_{j\in I}$ 
since condition 1 of Definition 3.2 does not need to hold if $\varphi$ is replaced by $\varphi_j$ and 
we take $a_1,a_m\in {\mathcal A}_{j,j}$.}
\end{Remark}   

Let us assume now that for any natural $n$ we have an $n$-dimensional square array of 
self-adjoint variables $(X_{i,j}(n))$ in a unital *-algebra ${\mathcal A}(n)$ equipped with 
an array of states $(\varphi_{i,j}(n))$ defined by a family of states  
$(\varphi_{j}(n))_{1\leq j\leq n}$. Blocks are defined by partitioning the set $\{1,2,\ldots ,n\}$ into disjoint 
non-empty subsets,   
$$
[n]=N_1\cup N_2\cup \ldots \cup N_r
$$ 
where $r\in {\mathbb N}$, with dependence on $n$ supressed in the notation, 
whose sizes increase as $n\rightarrow \infty$ so that
$n_{j}/n\rightarrow d_{j}$, where $n_j$ is the cardinality of $N_j$ for any $j \in [r]$.
Then the numbers $d_j$ form a diagonal matrix $D$ of trace one called the {\it dimension matrix}.

Concerning the distributions of the considered arrays, we assume that
\begin{enumerate}
\item[(A1)]
$(X_{i,j}(n))$ is matricially free with respect to 
$(\varphi_{i,j}(n))$ for any $n\in {\mathbb N}$,
\item[(A2)]
the variables have zero expectations,
$$
\varphi_{i,j}(n)(X_{i,j}(n))=0
$$ 
for all $i,j \in [n]$ and $n\in {\mathbb N}$,
\item[(A3)]
their variances are block-identical and are of order $1/n$, namely
$$
\varphi_{i,j}(n)(X_{i,j}^{2}(n))=\frac{u_{p,q}}{n}
$$
for any $i\in N_p, j\in N_q$, where each $u_{p,q}$ is a non-negative real number,
\item[(A4)]
their moments are uniformly bounded, i.e. $\forall m\;\exists M_{m}\geq 0$ such that 
$$
|\varphi_{i,j}(n)(X_{i,j}^{m}(n))|\leq \frac{M_{m}}{n^{m/2}}
$$
for all $i,j\in [n]$ and $n\in {\mathbb N}$.
\end{enumerate}
In particular, if the distributions of the $X_{i,j}(n)$'s in the states $\varphi_{i,j}(n)$ 
are block-identical, assumptions (A3)-(A4) are satisfied. 

Following [5], where we studied limit distributions of
random pseudomatrices under $\psi(n)$, we consider now normalized partial traces
$$
\psi_{k}(n)=\frac{1}{n_k}\sum_{j\in N_k}\varphi_j(n)
$$
where $k\in [r]$ and $n\in {\mathbb N}$. Easy modifications of the proofs given in  
[5, Lemma 6.1] and [5, Theorem 6.1] lead to combinatorial formulas for the limit distributions 
under partial traces given above. 
\begin{Theorem} $[{\rm 5}]$
Under assumptions (A1)-(A4), the limit distributions of random pseudomatrices under partial traces
have the form 
$$
\lim_{n\rightarrow \infty}\psi_{k}(n)(S^{m}(n))=\sum_{\pi\in \mathcal{NC}_{m}^{2}}b_{k}(\pi)
$$
where $k \in [r]$, $m\in {\mathbb N}$ and $B=DU$, with $D$ being the dimension matrix.
Consequently, $\psi(n)(S^{m}(n))$ converges to $\sum_{\pi\in \mathcal{NC}_{m}^{2}}b(\pi)$
as $n\rightarrow \infty$, where $b(\pi)=\sum_{k=1}^{r}d_{k}b_{k}(\pi)$.
\end{Theorem}
The above result refers to the `tracial' framework which reminds the limit theorem for random matrices
and is of main interest to us in this work. However, we have also shown in [5] that a similar 
result holds for the `standard' framework which reminds the central limit theorem
and involves the distributions of random pseudomatrices in the distinguished 
states $\varphi(n)$. This can be phrased as follows.
\begin{Theorem} $[{\rm 5}]$
Under assumptions (A1)-(A4) for arrays of states $(\varphi_{i,j}(n))$ 
defined by a distinguished state $\varphi(n)$ and the associated 
conditions $(\varphi_j(n))_{1\leq j\leq n}$ for each $n$, 
$$
\lim_{n\rightarrow \infty}\varphi(n)(S^{m}(n))=\sum_{\pi\in \mathcal{NC}_{m}^{2}}b_{0}(\pi)
$$
where $m\in {\mathbb N}$ and $B=DU$.
\end{Theorem}
In the sequel we will derive Hilbert space realizations of the limit 
joint distributions of random pseudomatrices $S(n)$ and their blocks 
under $\varphi(n)$, $\psi_k(n)$ and $\psi(n)$.
Interestingly enough, all these realizations are given 
on the same type of Hilbert space which is a 
matricially free product of Fock spaces.

\section{Matricially free Gaussian operators}

In this Section we shall introduce self-adjoint operators which play the role of 
matricially free Gaussian operators living in the matricially free product 
of Fock spaces. 
  
Recall that by the boolean and free Fock spaces over the Hilbert space ${\mathcal H}$, respectively, we understand
the direct sums 
$$
{\mathcal F}_{0}({\mathcal H})={\mathbb C}\xi \oplus {\mathcal H}\;\;\;{\rm and}\;\;\;
{\mathcal F}({\mathcal H})={\mathbb C}\xi \oplus \bigoplus_{m=1}^{\infty}{\mathcal H}^{\otimes m},
$$
where $\xi$ is a unit vector, endowed with the canonical inner products.
We shall use them to define an array of Fock spaces and their matricially free product.

\begin{Definition}
{\rm Let $({\mathcal H}_{i,j},\xi_{i,j})$ be an array of Hilbert spaces
with distinguished unit vectors. By the {\it matricially free product 
of $({\mathcal H}_{i,j},\xi_{i,j})$} we understand the pair $({\mathcal H}, \xi)$, where 
$$
{\mathcal H}= {\mathbb C}\xi \oplus\bigoplus_{m=1}^{\infty}
\bigoplus_{(i_{1},i_2)\neq \ldots \neq (i_{m},i_m)}
{\mathcal H}_{i_1,i_2}^{0}
\otimes 
{\mathcal H}_{i_2,i_3}^{0}
\otimes \ldots \otimes 
{\mathcal H}_{i_{m},i_{m}}^{0},
$$
with ${\mathcal H}_{i,j}^{0}={\mathcal H}_{i,j}\ominus {\mathbb C} \xi_{i,j}$ and 
$\xi$ being a unit vector, with the canonical inner product. We denote it
$({\mathcal H}, \xi)=*^{M}_{i,j}({\mathcal H}_{i,j},\xi_{i,j})$.}
\end{Definition}

We already know that the matricially free Fock space, 
which is a matricially free analog of the free Fock space, is the matricially free
product of an array of free Fock spaces [5]. Nevertheless, in order to 
find Hilbert space realizations of the limit laws
of Theorems 3.1 and 3.2, it suffices to take a matricially free product 
of an array of Fock spaces, in which free Fock spaces are put only on the diagonal 
and the boolean Fock spaces elsewhere. Clearly, we obtain a truncation of the 
matricially free Fock space in this fashion.
This structure is related to the difference between diagonal and off-diagonal random variables
and is related to the difference between multiplication of diagonal and off-diagonal 
blocks of usual matrices. 

\begin{Definition}
{\rm By the {\it matricially free-boolean Fock space} over the array 
$\widehat{\mathcal H}=({\mathcal H}_{i,j})$ we shall understand the matricially free product 
$$
({\mathcal F}, \xi)=*^{M}_{i,j}({\mathcal F}_{i,j}, \xi_{i,j}),\;\;\;
{\rm where} \;\;\;
{\mathcal F}_{i,j}=
\left\{
\begin{array}{ll}
{\mathcal F}({\mathcal H}_{j,j})& {\rm if}\;\;i=j\\
{\mathcal F}_{0}({\mathcal H}_{i,j})& {\rm if}\;\;i\neq j
\end{array}
\right. 
$$ 
and $\xi_{i,j}$ denotes the distinguished unit vector in ${\mathcal F}_{i,j}$.}
\end{Definition}

\begin{Remark}
{\rm 
If we have a square array 
$\widehat{\mathcal H}=({\mathcal H}_{i,j})$, where 
${\mathcal H}_{i,j}\cong {\mathcal H}_{i}$ for any 
$i,j\in I$ and $({\mathcal H}_{i})$ is a family of Hilbert spaces, 
then we have a natural isomorphism
$$
{\mathcal F}\cong {\mathcal F}(\bigoplus_{j}{\mathcal H}_{j})
$$
since each tensor product  
$$
({\mathcal H}_{i_1,i_1}^{\otimes (n_1-1)}\otimes {\mathcal H}_{i_1,i_2}\otimes 
{\mathcal H}_{i_2,i_2}^{\otimes (n_2-1)})\otimes \ldots 
\otimes 
({\mathcal H}_{i_{m-1},i_{m-1}}^{\otimes (n_{m-1}-1)}\otimes {\mathcal H}_{i_{m-1},i_{m}}\otimes 
{\mathcal H}_{i_m,i_m}^{\otimes n_m})
$$
is isomorphic to
$$
{\mathcal H}_{i_1}^{\otimes n_1}\otimes 
{\mathcal H}_{i_2}^{\otimes n_2}\ldots \otimes {\mathcal H}_{i_m}^{\otimes n_m}
$$
for any $i_1\neq i_2\neq \ldots \neq i_m$ and $m,n_1,n_2,\ldots , n_m\in {\mathbb N}$.
Thus, ${\mathcal F}$ is in this case also isomorphic to the strongly matricially free Fock space 
${\mathcal R}(\widehat{\mathcal H})$ introduced in [5].}
\end{Remark}

\begin{Example}
{\rm A simple example of a matricially free-boolean Fock space 
is the matricially free product of the two-dimensional array of the form
$$
({\mathcal F}_{i,j})
=
\left(
\begin{array}{ll}
l_{2}(G_{1,1}) & l_{2}(G_{1,2})\\
l_{2}(G_{2,1}) & l_{2}(G_{2,2})
\end{array}
\right),
$$
where $G_{j,j}=FS(g_{j,j})$, the free semigroup on one generator $g_{j,j}$, where $j\in \{1,2\}$,
and $G_{i,j}={\mathbb Z}_{2}$ with the generator denoted $g_{i,j}$ for $i\neq j$.
Then ${\mathcal F}=l_{2}(G)$, where $G$ is the `matricially free product of semigroups' 
$G_{i,j}$, by which we understand the subset of their free 
product $*_{i,j}G_{i,j}$ given by the union
$$
G=\bigcup_{n=0}^{\infty}G^{(n)}
$$
of disjoint subsets, where $G^{(n)}$ consists of words of type $g_1g_2\ldots g_n$ with 
$g_{k}$ being an element of $G_{i_k,j_k}^{0}$, where $G_{i,j}^{0}:=G_{i,j}\setminus \{\epsilon_{i,j}\}$, 
where $\epsilon_{i,j}$ is the unit in $G_{i,j}$, with 
$((i_1,j_1), \ldots , (i_n,j_n))\in \Lambda$ and $i_n=j_n$. For instance,
\begin{eqnarray*}
G^{(0)}&=&\{e\},\\
G^{(1)}&=&\{g_{1,1}^{k},\, g_{2,2}^{k}:\; k\in {\mathbb N}\},\\
G^{(2)}&=&\{g_{2,1}^{}g_{1,1}^{k},\, g_{1,2}^{}g_{2,2}^{k}:\; k\in {\mathbb N}\},\\ 
G^{(3)}&=&\{g_{2,2}^{k}g_{2,1}^{}g_{1,1}^{m},\, g_{1,1}^{k}g_{1,2}^{}g_{2,2}^{m},
\,g_{1,2}^{}g_{2,1}^{}g_{1,1}^{m}, \,g_{2,1}^{}g_{1,2}^{}g_{2,2}^{m}:\; k,m\in {\mathbb N} \},
\end{eqnarray*}
etc., with the remaining subsets consisting of words
built from `matricially free products' of powers of the generators,
where the diagonal generators admit all natural powers, 
whereas the off-diagonal ones admit only the powers equal to one. }
\end{Example}

\begin{Definition}
{\rm Let $A=(\alpha_{i,j})$ be a diagonal-containing array of positive real numbers
and let $({\mathcal H}_{i,j})=({\mathbb C}e_{i,j})$ be the associated array of
Hilbert spaces. By the {\it matricially free creation operators} associated with $A$ we understand
operators of the form
$$
\varsigma_{i,j}=\alpha_{i,j}\tau^{*}\ell(e_{i,j})\tau, 
$$
where $\tau: {\mathcal F}\rightarrow {\mathcal F}(\bigoplus_{i,j}{\mathcal H}_{i,j})$
is the canonical embedding and
the $\ell(e_{i,j})$'s denote the canonical free creation operators.
By the {\it matricially free annihilation operators} 
and the {\it matricially free Gaussian operators} we understand their adjoints $\varsigma_{i,j}^{*}$
and sums denoted $\zeta_{i,j}=\varsigma_{i,j}+\varsigma^{*}_{i,j}$, respectively. }
\end{Definition}

We shall assume now that $A$ is a diagonal-containing subarray of a finite square array
and that it is indexed by the set $J$ (thus, $\Delta\subseteq J \subseteq I\times I$).  
Moreover, it is convenient to assume that $A$ is the square root of another array $B$ 
taken entrywise, i.e. 
$$
\alpha_{i,j}=\sqrt{b_{i,j}}\;\;\;{\rm written}\;\;\;A=\sqrt{B}
$$
where $B$ can be assumed to be a subarray of a square array.
The dependence on $A$ of the operators defined above is supressed in our notations.

We shall prove below that the *-algebras ${\mathcal A}_{i,j}$, each generated by
$\varsigma_{i,j}$ and suitably defined unit $1_{i,j}$, respectively, are 
matricially free with respect to a suitably defined array of states
on $B({\mathcal F})$. 
Namely, $1_{i,j}$ is defined as the projection onto the subspace of 
${\mathcal F}$ onto which the *-algebra generated by the creation operator 
$\varsigma_{i,j}$ acts non-trivially. To be more precise, let us introduce 
projections $s_{i,j}$ and $r_{i,j}$ for any $(i,j)\in J$ which give an orthogonal decomposition of $1_{i,j}$, i.e.
$1_{i,j}=r_{i,j}+s_{i,j}$ and $r_{i,j}s_{i,j}=1$. 
Namely, $s_{i,j}$ is the canonical projection onto the subspace of ${\mathcal F}$ spanned by tensors which begin 
with $e_{i,j}$ for any $(i,j)\in J$. In turn, $r_{i,j}$ is the canonical projection onto
the subspace of ${\mathcal F}$ spanned by tensors which begin with $e_{j,k}$ for some $k$ 
if $i\neq j$, whereas $r_{j,j}$ is the canonical projection onto the subspace spanned by 
the vacuum vector and tensors which begin with $e_{j,k}$ for $k\neq j$, where $(j,k)\in J$.

Two types of transformations on the considered operators 
will be performed: truncations and symmetrizations.
Here, we shall introduce {\it truncated matricially free creation operators}
and {\it truncated units} by 
$$
\wp_{i,j}=\varsigma_{i,j}P \;\;\;{\rm and}\;\;\;
t_{i,j}=1_{i,j}P
$$
respectively, where $P$ is the canonical projection onto ${\mathcal F}\ominus {\mathbb C}\Omega$.
The {\it truncated matricially free annihilation operators} and 
{\it truncated matricially free Gaussian operators}, respectively, will be denoted
by $\wp_{i,j}^{*}$ and $\omega_{i,j}=\wp_{i,j}+\wp_{i,j}^{*}$ for any $i,j$.
For uniformity, one can put $P$ on both sides of the creation, annihilation and unit operators. 

Finally, in the case of finite dimensional arrays, which will be considered from now on,
we use the same symbol to denote the sum of operators in a given array, like
$$
\zeta=\sum_{(i,j)\in J}\zeta_{i,j}\;\;\;{\rm and}\;\;\;\omega=\sum_{(i,j)\in J}\omega_{i,j},
$$
called the {\it Gaussian pseudomatrix} and the {\it truncated Gaussian pseudomatrix}, respectively.
All adjoints are denoted in the usual way. 
Note that all these operators are bounded since all considered sums are finite.

\begin{Proposition}
Let $(\varphi_{i,j})$ be the array of states on $B({\mathcal F})$ defined 
by the vacuum state $\varphi$ and the family $(\psi_j)_{j\in I}$, 
where $\psi_j$ is the vector state associated with $e_{j,j}$ for any $j\in I$. Then
\begin{enumerate}
\item
the $\varphi$-distribution of $\zeta_{j,j}$ is the semicircle law of
radius $2\alpha_{j,j}$ for any $j$,
\item
the $\psi_j$-distribution of $\zeta_{i,j}$ is the Bernoulli law concentrated 
at $\pm \sqrt{2}\alpha_{i,j}$ for any $i\neq j$,
\item
the array $(\zeta_{i,j})$ is matricially free with respect to $(\varphi_{i,j})$.
\end{enumerate}
\end{Proposition}
{\it Proof.}
The first two claims follow easily from the definitions of the operators involved.
Next, instead of proving the third claim, we will prove a slightly more general result that
the array $({\mathcal A}_{i,j})$, where
${\mathcal A}_{i,j}={\mathbb C}\langle \varsigma_{i,j}, \varsigma_{i,j}^{*},1_{i,j}\rangle$
for any $(i,j)\in J$, is matricially free with respect to 
$(\varphi_{i,j})$.
Essentially, we proceed as in the free case [9], 
but we have slightly more complicated relations between the operators involved
(see also Proposition 4.2 of [5]).
In particular, one has to treat diagonal and off-diagonal subalgebras separately.
Let us consider first the off-diagonal case, when the creation and annihilation operators
satisfy relations
$$
\varsigma_{i,j}^{*}\varsigma_{i,j}=b_{i,j}r_{i,j}\;\;\;{\rm and }\;\;\;
\varsigma_{i,j}^{2}=0, \;\;\varsigma_{i,j}^{*2}=0
$$
for any $i\neq j$. In that case we have additional relations 
$$
r_{i,j}\varsigma_{i,j}=0, \;\;\;\;\varsigma_{i,j}r_{i,j}=\varsigma_{i,j}, \;\;\;\;
s_{i,j}\varsigma_{i,j}=\varsigma_{i,j},\;\;\;\;\varsigma_{i,j}s_{i,j}=0.
$$
Using all these relations and their adjoints, as well as equations
$$
\psi_{j}(r_{i,j})=1,\;\;\;\psi_{j}(s_{i,j})=0
$$
we deduce that an arbitrary noncommutative polynomial from 
${\mathcal A}_{i,j}\cap {\rm Ker}(\psi_{j})$, where $i\neq j$, is spanned by 
$\varsigma_{i,j}, \;\;\varsigma_{i,j}^{*}, \;\;{\rm and}\;\; s_{i,j}$.
In the diagonal case, the situation is similar to that in the free case since
each pair of diagonal creation and annihilation operators satisfies 
the relation
$$
\varsigma_{j,j}^{*}\varsigma_{j,j}=b_{j,j}1_{j,j}\;\;\;
{\rm with}\;\;\;\varphi(1_{j,j})=1,
$$
and therefore any polynomial from ${\mathcal A}_{j,j}\cap {\rm Ker}(\varphi)$ is spanned by 
$1_{j,j}$ and $\varsigma_{j,j}^{p}\varsigma_{j,j}^{*q}$, where $p+q>0$, for any $j$. 
Moreover, $s_{i,j}\varsigma_{k,l}=0$ and $\varsigma_{k,l}s_{i,j}=\delta_{l,i}\varsigma_{k,l}$
for any $(i,j)\neq (k,l)$. Hence, to prove matricial freeness of the array $({\mathcal A}_{i,j})$ 
with respect to $(\varphi_{i,j})$, it suffices to
show that 
$$
\varphi(\varsigma_{i_1,j_1}^{p_1}\varsigma_{i_1,j_1}^{*q_1}\ldots \varsigma_{i_m,j_m}^{p_m}\varsigma_{i_m,j_m}^{*q_m})=0
$$
for suitable powers $p_1,q_1, \ldots , p_m,q_m$ that depend on whether the corresponding ope\-rators are diagonal or not,
since $\varphi(s_{i_1,j_1}\ldots s_{i_m,j_m})=0$ for any off-diagonal pairs $(i_1,j_1)\neq \ldots \neq (i_m,j_m)$.
At this point we can use the same inductive argument as in the free case [9], which
gives the above `freeness condition'. Moreover, the definition of each $1_{i,j}$ shows that 
it is the projection onto the subspace onto which the *-algebra generated by $\varsigma_{i,j}$ 
acts non-trivially, which implies that the array $(1_{i,j})$ is the matricially free array of units.
\hfill $\blacksquare$\\

Note that in Proposition 4.1 each $\psi_j$ is a condition associated with $\varphi$ since there exists 
$b_{j}\in {\mathcal A}_{j,j}\cap{\rm Ker}\varphi$, namely $b_j=\varsigma_{j,j}/\alpha_{j,j}$, such that
$\psi_{j}(a)=\varphi(b_{j}^{*}ab_{j})$ for any $a\in B({\mathcal F})$. 
However, a similar result is obtained for arrays of states defined by the family 
$(\psi_j)_{j\in I}$, where $\psi_j$ is now the vector state 
on $B({\mathcal F}\ominus {\mathbb C}\Omega)$ associated with $e_{j,j}$
for any $j\in I$.

\begin{Proposition}
Let $(\psi_{i,j})$ be the array of states on $B({\mathcal F}\ominus {\mathbb C}\Omega)$ 
defined by the family $(\psi_j)_{j \in I}$, where
$\psi_j$ is the vector state associated with $e_{j,j}$ for any $j$. Then
\begin{enumerate}
\item
the $\psi_j$-distribution of $\omega_{j,j}$ is the semicircle law of
radius $2\alpha_{j,j}$ for any $j$,
\item
the $\psi_j$-distribution of $\omega_{i,j}$ is the Bernoulli law concentrated 
at $\pm \sqrt{2}\alpha_{i,j}$ for any $i\neq j$,
\item
the array $(\omega_{i,j})$ is matricially free with respect to 
$(\psi_{i,j})$.
\end{enumerate}
\end{Proposition}
{\it Proof.}
The proof is similar to that of Proposition 4.1.
\hfill $\blacksquare$\\

\section{Fock-space realizations of limit distributions}

Using the matricially free Gaussian operators and their truncations, 
we will find realizations on ${\mathcal F}$ of the limit distributions of Theorems 3.1-3.2.

To each $\pi\in {\mathcal NC}_{m}^{2}$ for $m$ even 
we assign natural products of creation and annihilation operators.
First, to each $\pi\in \mathcal{NC}_{m}^{2}$ we assign the sequence $\epsilon(\pi)=(\epsilon_{1},\epsilon_{2}, \ldots , \epsilon_{m})$, where
$\epsilon_{j}=1$ whenever $j\in {\mathcal R}(\pi)$ and $\epsilon_{k}=*$
whenever $k\in {\mathcal L}(\pi)$. Then, 
to each $\pi\in \mathcal{NC}_{m}^{2}$ we assign
products of creation and annihilation operators
$$
\varsigma(\pi)=\varsigma^{\epsilon_{1}}\varsigma^{\epsilon_{2}}\ldots \varsigma^{\epsilon_{m}}\;\;\;
{\rm and}\;\;\;
\wp_{j}(\pi)=\wp^{\epsilon_{1}}\wp^{\epsilon_{2}}\ldots \wp^{\epsilon_{m}}
$$
for any $j\in [r]$,
where $(\epsilon_1,\epsilon_2, \ldots , \epsilon_m)=\epsilon(\pi)$.

There is a nice relation between the expectations of these products 
and numbers $b_{j}(\pi)$ defined in Section 2. 
In the case of $b_{0}(\pi)$ we will use 
the expectations in the vacuum state $\varphi$ and 
for the remaining $j$'s we take $\psi_{j}$'s associated with vectors $e_{j,j}$.
Finally, to obtain $b(\pi)$, we will use the convex linear combination of states of the form 
$$
\psi=\sum_{j=1}^{r}d_{j}\psi_{j}
$$
where numbers $d_1,d_2,\ldots , d_r$ are taken from the dimension matrix $D$.

We are ready to give Fock-space realizations of the limit distributions of random pseudomatrices 
in terms of the distributions of Gaussian and truncated Gaussian pseudomatrices.

\begin{Lemma}
For the limit distributions of Theorem 3.2, it holds that
$$
\sum_{\pi\in \mathcal{NC}_{m}^{2}}b_{0}(\pi)=\varphi
(\zeta^{m})\;\;\;
$$
for any $m\in {\mathbb N}$, where $\zeta$
is the Gaussian pseudomatrix.
\end{Lemma}
{\it Proof.}
If $m$ is odd, both sides of the first equation are clearly zero. 
Therefore, suppose that $m$ is even.
We have
$$
\varphi(\zeta^{m})
=
\sum_{\epsilon_1, \epsilon_2, \ldots , \epsilon_m\in \{1,*\}}
\varphi(\varsigma^{\epsilon_1}\varsigma^{\epsilon_2}\ldots \varsigma^{\epsilon_m})
=
\sum_{\pi\in \mathcal{NC}^{2}_{m}}\varphi(\varsigma(\pi))
$$
where we used the fact that products of creation and annihilation operators
corresponding to sequences $(\epsilon_1,\epsilon_2, \ldots , \epsilon_m)$ which are not 
associated with non-crossing pair partitions give zero contribution (this follows
from the definition of the $\varsigma_{i,j}^{\epsilon}$). 
The assertion of the lemma
follows then from the combinatorial formula 
$$
b_{0}(\pi)=\varphi(\varsigma(\pi))
$$
for each $\pi\in \mathcal{NC}_{m}^{2}$. To prove it, first observe that
$$
\varphi(\varsigma_{p_1,q_1}^{\epsilon_{1}}\ldots \varsigma_{p_m,q_m}^{\epsilon_m})
=0 \;\;\;{\rm unless}\;\;(\epsilon_1, \ldots , \epsilon_m)=\epsilon(\pi)
$$
for some $\pi\in \mathcal{NC}_{m}^{2}$. Now, if $\pi\in \mathcal{NC}_{m}^{2}(p_1,q_1, \ldots , p_m,q_m)$
for some $p_1,q_1, \ldots , p_m,q_m\in [r]$ and $(\epsilon_1, \ldots , \epsilon_m)=\epsilon(\pi)$, then
we claim that
$$
\varphi(\varsigma_{p_1,q_1}^{\epsilon_{1}}\ldots \varsigma_{p_m,q_m}^{\epsilon_m})
=b_{0}(\pi,f)
$$
where $f$ is the coloring of $\pi$ defined by $p_k$, $k\in \mathcal{L}(\pi)$.
This claim can be proved by induction. 
Suppose that the last annihilation operator in this mixed moment
is indexed by $k$, i.e. $\epsilon_{k}=*$. 
Then $\{k,k+1\}$ must be a block which does not have any inner blocks.
The corresponding product of operators is of the form 
$\varsigma_{i_k,j_k}^{*}\varsigma_{i_k,j_k}^{}$.
Two cases are possible:
\begin{enumerate}
\item
If $k=m-1$ and $p_k=q_k$, then 
$\varsigma_{p_k,q_k}^{*}\varsigma_{p_k,q_k}^{}$ acts on $\Omega$ and gives 
$b_{p_k,q_k}\Omega$.
\item
If $k<m-1$ and $q_{k+1}=p_{k+2}$, then 
$\varsigma_{p_k,q_k}^{*}\varsigma_{p_k,q_k}^{}$ acts on some simple tensor $h$ and gives $b_{p_k,q_k}h$.
Here, $q_k$ colors the nearest outer block of $\{k,k+1\}$.
\end{enumerate}
If we repeat this procedure
for the product of operators corresponding to the 
partition $\pi'$ obtained from $\pi$ by removing block $\{k,k+1\}$, 
we obtain the product of $b_{p,q}$'s which appears in the combinatorial formula for $b_{0}(\pi,f)$
after a finite number of steps. 
If we fix $\pi$ and sum over all $p_1,q_1, \ldots , p_m,q_m$ to which
$\pi$ is adapted (Definition 2.2), we obtain in fact the summation over $p_k, k\in \mathcal{L}(\pi)$,
equivalent to the summation over $F_{r}(\pi)$, which proves the formula for $b_{0}(\pi)$
since mixed moments corresponding to those $p_1,q_1,\ldots, p_m,q_m$ 
to which $\pi$ is not adapted are equal to zero.
This completes the proof.
\hfill $\blacksquare$

\begin{Lemma}
For the limit distributions of Theorem 3.1, it holds that
$$
\sum_{\pi\in \mathcal{NC}_{m}^{2}}b_{k}(\pi)=\psi_{k}
(\omega^{m})
$$
for $k \in [r]$ and $m\in {\mathbb N}$
and hence $\sum_{\pi\in \mathcal{NC}_{m}^{2}}b(\pi)=\psi(\omega^{m})$,
where $\omega$ is the truncated Gaussian pseudomatrix.  
\end{Lemma}
{\it Proof.}
The proof is similar to that of Lemma 5.1 and is based on the analogous 
combinatorial formula
$$
b_{k}(\pi)=\psi_{k}(\wp(\pi))
$$
for any $1 \leq k \leq r$ and $\pi\in \mathcal{NC}_{m}^{2}$, where $m$ is even and positive. 
The proof of that formula reduces to showing that
if $(\epsilon_1, \ldots , \epsilon_m)=\epsilon(\pi)$
for some $\pi\in \mathcal{NC}_{m}^{2}(p_1,q_1, \ldots , p_m,q_m)$ and $k=q_m$, then
$$
\psi_{k}(\wp_{p_1,q_1}^{\epsilon_{1}}\ldots \wp_{p_m,q_m}^{\epsilon_m})
=
b_{k}(\pi,f)
$$
where $f$ is the coloring of $\pi$ defined by $p_k$, $k\in {\mathcal L}(\pi)$, 
with $k$ coloring the imaginary block. 
As compared with the proof for $\varphi$, instead of acting on $\Omega$, we need to act
on $e_{k,k}^{}$. However, thanks to the projection $P$ 
onto ${\mathcal F}\ominus {\mathbb C}\Omega$, which appears in the definition
of the $\wp_{k,k}^{\epsilon}$, each $e_{k,k}^{}$
plays the role of the vacuum vector with respect to the $\wp_{i,k}^{\epsilon}$
for any $i,k$, including the case when $i=k$, since $\wp_{k,k}^{*}e_{k,k}^{}=0$. 
Therefore, the arguments are similar to those for $\varphi$, 
except that $\wp_{i,k}e_{k,k}$ is non-zero
for arbitrary $i$. In terms of diagrams, this means that each covering block 
of $\pi$ contributes $b_{i,k}$ for various $i$'s. 
Finally, when we use the definition of $\psi$ and
sum over $k\in [r]$, we obtain for each $k$ the extra factor
$d_{k}$ which appears in the combinatorial formula for $b(\pi)$.
Consequently, $b(\pi)=\psi(\wp(\pi))$ for any such $\pi$, which then leads to the 
formula for $\psi(\omega^n)$.
\hfill $\blacksquare$

\begin{Example}
{\rm Consider the moment corresponding to the partition $\pi\in \mathcal{NC}_{4}^{2}$ 
consisting of two blocks: $\{1,4\}$ and $\{2,3\}$. Let $A=\sqrt{B}\in M_{2}({\mathbb R})$ be
the square root taken entrywise, where $B=DU$ and we set $\alpha_{1,1}=\alpha, \alpha_{1,2}=\beta, \alpha_{2,1}=\gamma, 
\alpha_{2,2}=\delta$. Then
$$
b_{0}(\pi)
=
\alpha^{4}+ \alpha^{2}\gamma^{2}+\delta^{2}\beta^{2}+\delta^{4}.
$$
On the other hand, 
\begin{eqnarray*}
\varsigma^{2}\Omega&=&\alpha^2(e_{1,1}\otimes e_{1,1}) + \alpha \gamma (e_{2,1}\otimes e_{1,1})\\
&+&
\beta \delta (e_{1,2}\otimes e_{2,2}) + \delta^2(e_{2,2}\otimes e_{2,2})
\end{eqnarray*}
and thus the `symmetric' action of $(\varsigma^{*})^{2}$ gives exactly $b_{0}(\pi)\Omega$.}
\end{Example}
\begin{Example}
{\rm For the same partition $\pi$ as in the above example we obtain 

\begin{eqnarray*}
b(\pi)
&=&
d_1
(\alpha^4
+
\alpha^2\gamma^2 + \beta^2\gamma^2 + \gamma^2\delta^2)\\
&+&
d_2(\alpha^2\beta^2+ \beta^2\gamma^2 +\beta^2\delta^2
+
\delta^4)
\end{eqnarray*}
On the other hand,
$$
\psi(\wp(\pi))=d_1\psi_{1}(\wp(\pi))+ d_2\psi_{2}(\wp(\pi))
$$
where $\wp(\pi)=\wp^*\wp^*\wp\wp$, and we get
\begin{eqnarray*}
\wp^{2}e_{1,1}
&=&
\gamma
\left(\delta(e_{2,2}\otimes e_{2,1}\otimes e_{1,1})
+ 
\beta(e_{1,2}\otimes e_{2,1}\otimes e_{1,1})\right)\\
&+&
\alpha
\left(\alpha(e_{1,1}\otimes e_{1,1}\otimes e_{1,1})
+
\gamma(e_{2,1}\otimes e_{1,1}\otimes e_{1,1})\right)
\end{eqnarray*}
which, in view of the `symmetric' action of the adjoint, gives
$$
d_1\psi_{1}(\wp^{*}\wp^{*}\wp \wp)=
d_1(\gamma^2\delta^2 + \beta^2\gamma^2 + \alpha^2\gamma^2 + \alpha^4).
$$
A similar expression is obtained for $d_2\psi_{2}(\wp^{*}\wp^{*}\wp \wp)$ 
with $\alpha$ interchanged with $\delta$ and $\beta$ interchanged with $\gamma$. 
The sum of both expressions agrees with $b(\pi)$. }
\end{Example}

\section{Matricial freeness of blocks}

We show in this section that blocks of finite-dimensional arrays 
of matricially free random variables are matricially free with respect to
an appropriately defined array of states. This is the analog of the property 
of free random variables which says that families of sums of free random 
variables are free. 

\begin{Definition}
{\rm Suppose that the $n$-dimensional array $(X_{i,j})$ of variables from a unital algebra ${\mathcal A}$ 
is matricially free with respect to $(\varphi_{i,j})$ and let 
$[n]=N_{1}\cup N_{2}\cup  \ldots \cup N_{r}$ be a partition of $[n]$ ino disjoint non-empty 
subsets. The sums of the form 
$$
S_{p,q}=\sum_{(i,j)\in N_p \times N_q}X_{i,j}
$$
where $p,q \in [r]$, will be called {\it blocks} of the pseudomatrix $S=\sum_{i,j}X_{i,j}$.
We will say that the array $(X_{i,j})$ has {\it block-identical distributions} 
with respect to $(\varphi_{i,j})$ if the $\varphi_{i,j}$-distributions
of $X_{i,j}$ are the same for all $(i,j)\in N_p\times N_q$ and fixed $p,q\in [r]$.}
\end{Definition}

We need to define the associated array of block units which 
satisfy the conditions of Definition 3.1. In order to construct them
on the level of noncommutative probability spaces, 
let us first return to the more intuitive framework of the matricially free 
product of representations of the array $({\mathcal A}_{i,j})$ of 
unital $C^{*}$-algebras on the matricially free product of Hilbert spaces 
$$
({\mathcal H},\xi)=*_{i,j}^{M}({\mathcal H}_{i,j},\xi_{i,j})
$$ 
and then carry over the corresponding definition to 
the algebraic framework of noncommutative probability spaces.
We assume that each ${\mathcal A}_{i,j}$ is equipped with an internal unit $1_{i,j}$ and a state $\varphi_{i,j}$.
It is the matricially free product of Hilbert spaces on which one defines
the canonical *-representations of $({\mathcal A}_{i,j}, \varphi_{i,j})$. Namely, let  
$({\mathcal H}_{i,j},\pi_{i,j},\xi_{i,j})$ be the associated GNS triples, so that
$\varphi_{i,j}(a)=\langle \pi_{i,j}(a)\xi_{i,j},\xi_{i,j}\rangle$ for any $a\in {\mathcal A}_{i,j}$.
By $\lambda_{i,j}$ we denote the canonical *-representation of $({\mathcal A}_{i,j}, \varphi_{i,j})$ 
on $({\mathcal H},\xi)$. Using these representations and appropriate partial isometries,
we have defined in [5] their matricially free product $\lambda=*_{i,j}^{M}\pi_{i,j}$
which maps $\sqcup_{i,j}{\mathcal A}_{i,j}$ into $B({\mathcal H})$.

In this connection note that $\lambda$ does not send the units $1_{i,j}$ onto the unit of $B({\mathcal H})$. In fact, 
$\lambda(1_{i,j})=r_{i,j}+s_{i,j}$, where 
\begin{eqnarray*}
r_{i,j}&=&{\rm projection}\;{\rm onto}\;\;{\mathcal H}(i,j)\\
s_{i,j}&=&{\rm projection}\;{\rm onto}\;\;{\mathcal K}(i,j)
\end{eqnarray*}
where 
${\mathcal K}(i,j)=\mathcal{H}_{i,j}^{0}\otimes \mathcal{H}(i,j)$
and 
${\mathcal H}(i,j)$ denotes the subspace of ${\mathcal H}$ onto which the left 
free action of $\lambda({\mathcal A}_{i,j})$ is non-trivial. By the left free action of $\lambda({\mathcal A}_{i,j})$ we understand the left action onto the subspace 
spanned by simple tensors which do not begin with vectors from ${\mathcal H}_{i,j}^{0}$ and, in addition, by $\xi$ if $i=j$.
Clearly, $r_{i,j}\perp s_{i,j}$ for any fixed $i,j$ and 
thus their sum is the canonical projection onto the subspace of 
${\mathcal H}$ onto which $\lambda({\mathcal A}_{i,j})$ acts non-trivially.

In the propositions given below we state basic properties of the projections involved that 
are needed for the construction of block units. Note that the case of ${\rm card}(I)=1$ is trivial from the point of view 
of matricially free product structures and that is why it is not treated. 
By ${\mathbb C}[a_i, i\in I]$ we denote the unital algebra of polynomials in the commuting 
indeterminates $a_{i}, i\in I$. Recall that ${\mathcal I}$ stands for the commutative unital algebra 
generated by the units $\lambda(1_{i,j})\equiv 1_{i,j}$. 

\begin{Proposition}
If $\;{\rm card}(I)>1$, then ${\mathcal I}={\mathbb C}[r_{i,j},s_{i,j},p_{\xi}: i,j\in I]$,
where $p_{\xi}$ is the canonical projection onto ${{\mathbb C}\xi}$. 
\end{Proposition}
{\it Proof.}
Observe that we have relations

$$
1_{i,i}1_{j,j}=p_{\xi},\;\;\;
1_{j,j}1_{j,k}=s_{j,k}\;\;\;{\rm and}\;\;\;
1_{j,j}1_{k,j}=s_{j,j}
$$
whenever $i\neq j \neq k$, which implies that $s_{j,k},p_{\xi} \in {\mathbb C}[1_{i,j}:i,j\in I]$
for any $j,k$, and since $r_{j,k}=1_{j,k}-s_{j,k}$, it holds that 
${\mathbb C}[r_{i,j},s_{i,j},p_{\xi}: i,j\in I] \subseteq {\mathbb C}[1_{i,j}:i,j\in I]$. 
The reverse inclusion is obviously true, hence the proof is completed.
\hfill $\blacksquare$

\begin{Proposition} 
If ${\rm card}(I)>1$ and $J_1,J_2\subseteq I$ are identical or disjoint, then the algebra ${\mathbb C}[1_{i,j}:i\in J_1,j\in J_2]$ contains the canonical projection ${\bf 1}_{J_1,J_2}$ onto the subspace of ${\mathcal H}$ onto which the algebra generated by $\{\lambda(a_{i,j}):a_{i,j}\in \mathcal A_{i,j},\;i\in J_1,j\in J_2\}$ 
acts non-trivially.
\end{Proposition}
{\it Proof.}
We want to construct a projection ${\bf 1}_{J_1,J_2}$ onto the subspace of ${\mathcal H}$ that would be 
suitable for the left action of the $\lambda(a_{i,j})$, where $i\in J_1$ and $j\in J_2$.
If $J_1\cap J_2=\emptyset$, then this subspace is the orthogonal direct sum
$$
{\mathcal H}_{J_1,J_2}=
\bigoplus_{i\in J_1,j\in J_2}{\mathcal K}(i,j)
\oplus
\bigoplus_{j\in J_2}\left({\mathcal K}(j,j)
\oplus 
{\mathcal H}(j,j)\ominus {\mathbb C}\xi\right)
$$
and the associated canonical projection is 
\begin{eqnarray*}
{\bf 1}_{J_1,J_2}&=&\sum_{i\in J_1,j\in J_2}s_{i,j}+\sum_{j\in J_2}(s_{j,j}+r_{j,j}-p_{\xi})\\
&=&
\sum_{i\in J_1,j\in J_2}(1_{i,j}-1_{j,j}+1_{i,i}1_{j,j})+
\sum_{j\in J_2}(1_{j,j}-1_{j,j}1_{k,k})
\end{eqnarray*}
where $k\in J_{1}$ is arbitrary. In turn, if $J_1= J_2=J$, then we obtain the orthogonal direct sum
$$
{\mathcal H}_{J,J}={\mathbb C}\xi
\oplus
\bigoplus_{j\in J}
\left(
{\mathcal K}(j,j)
\oplus {\mathcal H}(j,j)\ominus {\mathbb C}\xi\right)
$$
and the corresponding canonical projection takes the form
\begin{eqnarray*}
{\bf 1}_{J,J}&=&\sum_{j\in J}(s_{j,j}+r_{j,j}-p_{\xi})+p_{\xi}\\
&=&
\sum_{j\in J}(1_{j,j}-1_{j,j}1_{k,k}) + 1_{k,k}1_{l,l},
\end{eqnarray*}
where indices $k,l$ are such that $k\neq l$ and otherwise are arbitrary elements of $J$.
We have used the fact that $p_{\xi}=1_{i,i}1_{j,j}$ for any $i\neq j$.
This completes the proof.
\hfill $\blacksquare$\\

The proof of the above proposition enables us to construct block units which are internal units in 
the commutative algebras ${\mathbb C}[S_{p,q},{\bf 1}_{p,q}]$, where $p,q\in [r]$, 
each generated by a block and the corresponding block unit. 
Namely, we set 
$$
{\bf 1}_{p,q}:={\bf 1}_{N_p,N_q}
$$ 
for any $p,q \in [r]$, where the right hand side is defined in terms of $1_{i,j}$'s 
in exactly the same way as in the proof of Proposition 6.2.

It will also be useful to introduce the following terminology. Namely, if 
$$
a_k\in {\mathcal A}_{i_k,j_k}\cap {\rm Ker}(\varphi_{i_k,j_k})\;\; {\rm for}\;\; 1\leq k \leq n\;\;{\rm and}
\;\;((i_1,j_1), \ldots , (i_n,j_n))\in \Lambda ,
$$
we will say that the product $a_1a_2\ldots a_n$ is in the {\it matricially free kernel form} with respect to $(\varphi_{i,j})$.
\begin{Theorem}
Let $(\varphi_{i,j})$ be the array of states on ${\mathcal A}$ defined by 
the family $(\varphi_{j})_{1\leq j \leq n}$,  
and let $(\psi_{p,q})$ be defined by the family of associated 
normalized partial traces $(\psi_{q})_{1\leq q \leq r}$.
If $(X_{i,j})$ is matricially free (strongly matricially free) and has block-identical distributions 
with respect to $(\varphi_{i,j})$, then $(S_{p,q})$ is matricially free 
with respect to $(\psi_{p,q})$.
\end{Theorem}
{\it Proof.}
We will assume that $(X_{i,j})$ is matricially free 
since the proof for strong matricial freeness is analogous.
Clearly, the unital algebra generated by $({\bf 1}_{p,q})$ is commutative. 
We claim that the array $({\bf 1}_{p,q})$ is a matricially free array of units associated with
$(S_{p,q})$ and $(\psi_{p,q})$. The proof of condition (1) of Definition 3.1 for $(\psi_{p,q})$ 
follows easily from the same condition for $(\varphi_{i,j})$.
To prove condition (2), we need to evaluate mixed momeints of type 
$$
\varphi_j(w{\bf 1}_{p,q}w_1w_{2}\ldots w_m)
$$
where $w_{k}\in {\rm Ker}(\psi_{p_k,q_k})$ 
is a polynomial in $S_{p_k,q_k}$ for each $k\in [m]$ and $j\in [n]$
and the product $w_1w_2\ldots w_m$ is in the matricially free kernel form with respect
to $(\psi_{p,q})$.  
We shall reduce the computations to moments of type 
$\varphi_j(w{\bf 1}_{i,j}v_1v_2\ldots v_h)$, where $v_1v_2\ldots v_{h}$
is in the (strongly) matricially free kernel form with respect to $(\varphi_{i,j})$.
For that purpose, let us express each power of $S_{p_k,q_k}$ which appears in $w_{k}$
in terms of variables which are in the kernels of the $\varphi_{i,j}$. This procedure,
described in more detail below, is applied to $w_{m}, w_{m-1} \ldots , w_{1}$ (in that order).

Let $w(X)=\sum_{r=0}^{s}c_{r}X^r$ be an arbitrary polynomial and 
let $S_{p,q}=\sum_{(i,j)\in N_p\times N_q}X_{i,j}$.
We decompose each positive power of $X_{i,j}$ which appears in $w(S_{p,q})$ as
$$
X_{i,j}^{n}=(X_{i,j}^{n})^{0}+\varphi_{i,j}(X_{i,j}^{n})1_{i,j},
$$
where $i\in N_{p},j\in N_{q}$. Using the fact that $(1_{i,j})$ is a matricially 
free array of units, we observe that in the computations of the mixed moments of 
the given type we can use, without loss of generality, polynomials of the form
$$
w(S_{p,q})=\sum_{r=1}^{s}\sum_{i_1, \ldots , i_{r+1}}\sum_{n_1, \ldots, n_r}
c_{i_1, \ldots , i_r, i_{r+1}}^{n_1, \ldots, n_r}
(X_{i_1,i_2}^{(n_1)})^{0}
\ldots 
(X_{i_r,i_{r+1}}^{(n_r)})^{0}
+ 
c_{p,q}\sum_{(i,j)\in N_{p}\times N_{q}}1_{i,j}, 
$$
with summations over $i_1,\ldots, i_{r+1}$ and $n_1, \ldots , n_r$ 
run over some finite sets of natural numbers, with $n_1+\ldots +n_r\leq {\rm deg}(w)$,
since the remaining terms will give zero contribution to the considered moment
if the product of variables standing to the right of this polynomial
is in the matricially free kernel form with respect to $(\varphi_{i,j})$, 
which is the case if we carry out our computations going from the right to the left.

Note that if $p\neq q$, then $s=1$, but if $p=q$, then 
$s \leq {\rm deg}(w)$.
Let us also point out that thanks to our assumption
that the array $(X_{i,j})$ has block-identical distributions with respect to $(\varphi_{i,j})$, 
the same constant $c_{p,q}$ stands by each $1_{i,j}$ for any $i\in N_p,j\in N_q$, i.e. 
it does not depend on $i,j$. Moreover,
$$
\psi_{p,q}(\sum_{(i,j)\in N_p\times N_q}1_{i,j})=\frac{1}{n_{q}}\sum_{k\in N_{q}}\varphi_{k}(\sum_{(i,j)\in N_p \times N_q}1_{i,j})
=\frac{1}{n_{q}}\sum_{i,j}\varphi_{j}(1_{i,j})=n_{p}
$$
for any $p,q$, whereas the first sum in the above expression for $w(S_{p,q})$ belongs to 
${\rm Ker}(\psi_{p,q})$. 
These arguments lead us to the conclusion that in the computations 
of mixed moments of the considered type
we can take each polynomial $w_{k}$ to be of the form $w(S_{p_k,q_k})$ given above, with 
$c_{p_k,q_k}=0$ for each $k\in [m]$ since the product $w_1w_2\ldots w_m$ is assumed 
to be in the matricially free kernel form with respect to $(\varphi_{i,j})$.  

Therefore, each moment of type
$\varphi_j(w{\bf 1}_{p,q}w_1w_{2}\ldots w_m)$ is a sum of
moments of type $\varphi_j(w{\bf 1}_{p,q}v_1v_2\ldots v_h)$, where the product 
$v_1v_2\ldots v_h$ is in the matricially free kernel form
with respect to $(\varphi_{i,j})$. Similarly, each moment of type 
$\varphi_j(ww_1w_{2}\ldots w_m)$ is a corresponding sum of moments of type
$\varphi_j(wv_1v_2\ldots v_h)$ since the reduction described above does not depend on what
stands before the product $w_1w_2\ldots w_m$. 
It remains to observe that under each $\varphi_{j}$, the block unit 
${\bf 1}_{p,q}$ acts as the projection onto the linear span of products 
$v_1v_2\ldots v_j$ which are in the matricially 
free kernel form with respect to $(\varphi_{i,j})$ and begin
with $v_1\in {\mathcal A}_{i_1,j_1}$, where $i_1\in N_{q}$.
The proof of that fact follows from the definition of ${\bf 1}_{p,q}$ 
expressing it in terms of units $1_{i,k}$ which, under $\varphi_j$, 
act as projections onto the linear span of
$v_1v_2\ldots v_m$ which are in the matricially free kernel form
with respect to $(\varphi_{i,j})$ and begin with $v_1\in {\mathcal A}_{i_1,j_1}$,
where $i_1=k$. Of course, we use here the fact that $(p,q)\neq (p_1,q_1)$.
This completes the proof of condition (2) of Definition 3.1.

Moreover, using the normalization conditions for $1_{i,j}$'s, we obtain the normalization conditions
$$ 
\psi_r({\bf 1}_{p,q})=\frac{1}{n_r}\sum_{j\in N_r}\varphi_j({\bf 1}_{p,q})=\delta_{r,q}
$$
which completes the proof that $({\bf 1}_{p,q})$ is a matricially free array of units. 

Finally, the proof of condition (1) of Definition 3.2 is similar to that of condition (2) of Definition 3.1 presented above and is based on reducing the computations of the mixed moments 
$\varphi_j(w_1w_2\ldots w_n)$, where $w_1w_2\ldots w_n$ is in the matricially free kernel form 
with respect to $(\psi_{p,q})$, to mixed moments of type $\varphi_j(v_1v_2\ldots v_h)$, where
$v_1v_2\ldots v_h$ is in the matricially free kernel form with respect to $(\varphi_{i,j})$. This completes the proof.      
\hfill $\blacksquare$

\section{Asymptotic matricial freeness of blocks}

In this section we study the asymptotic joint distributions of blocks of random pseudomatrices.
In particular, we show that they are `asymptotically matricially free', which
is a notion analogous to asymptotic freeness. This generalizes the results of Section 6, where
we proved matricial freeness of blocks in the case when the array of matricially free variables 
has block-identical distributions.

For that purpose we will use a realization on the Fock space ${\mathcal F}$.
Having established the realization of the limit laws  
of random pseudomatrices $S(n)$ under normalized partial traces in Lemma 5.1, it is natural to
expect that it can be carried over to the level of {\it blocks} of $S(n)$ of the form 
$$
S_{p,q}(n)=\sum_{(i,j)\in N_p\times N_q}X_{i,j}(n)
$$ 
for any $n$, where we require that the arrays $(X_{i,j}(n))$
of self-adjoint random variables satisfy the assumptions of Theorem 3.1.

Before we proceed with examining the limit joint distribution of these blocks,
we define the notion of asymptotic matricial freeness, 
following the analogous notion of asymptotic freeness.
By the {\it block units} we shall understand units ${\bf 1}_{p,q}(n)$ defined for each $n$ in terms of
$1_{i,j}(n)$'s in exactly the same way as in Section 6.
In addition to the states $\varphi(n)$ and $\psi(n)$ on the algebras 
$\mathcal{A}(n)$, we will also use normalized partial traces 
$\psi_q(n)$, where $q \in [r]$ and $n\in {\mathbb N}$, defined as
before in the case of fixed $n$. 
Informally, the states $\psi_q(n)$ will play the role of 
conditions which converge to the conditions 
$\psi_{q}$ of Section 4 as $n\rightarrow \infty$. 

\begin{Definition}
{\rm Let $(\eta_{p,q}(n))$ be an $r$-dimensional array of functionals on 
the algebra of noncommutative polynomials 
${\mathbb C}\langle S_{p,q}(n), {\bf 1}_{p,q}(n), p,q \in [r]\rangle$.
We will say that $(S_{p,q}(n))$ is {\it asymptotically matricially free} 
with respect to $(\eta_{p,q}(n))$ if these functionals have pointwise limits
$(\eta_{p,q})$ as $n\rightarrow \infty$ 
with respect to which the limit array is matricially free.}
\end{Definition}

\begin{Theorem}
Under the assumptions of Theorem 3.1, the joint $\psi_{q}(n)$-distributions of 
blocks and block units converge 
to the joint $\psi_{q}$-distributions of the truncated matricially free Gaussian operators and  
the truncated units, respectively, as $n\rightarrow \infty$.
\end{Theorem}
{\it Proof.}
This result is a refinement of Theorem 3.1. The combinatorial arguments referring
to non-crossing partitions used in [5, Lemma 6.1] can be repeated except that since we now take products of the 
$S_{p,q}(n)$'s instead of powers of $S(n)$, all indices run only over some of the subsets $N_1,N_2, \ldots, N_r$ of 
the set $[n]$. Moreover, it suffices to consider the case of $m$ even, say $m=2s$, since
if $m$ is odd, both sides are zero, by standard arguments. 
We will first show that 
$$
\lim_{n\rightarrow \infty}
\psi_{q}(n)(S_{p_1,q_1}(n)\ldots S_{p_m,q_m}(n))
=
\psi_{q}(\omega_{p_1,q_1}\ldots \omega_{p_m,q_m})
$$
for any $q\neq 0$. 
Thus, for given $q$ and  
$p_1,q_1, \ldots , p_m,q_m$, in order to compute the left-hand side of the 
above equation, it suffices to take into account only the mixed moments
$$
\varphi_{j}(n)(X_{i_1,j_1}(n) \ldots X_{i_m,j_m}(n))
$$
in which $(i_k,j_k)\in N_{p_k}\times N_{q_k}$ for any $k\in [m]$ and $j\in N_q$.

However, in the limit $n\rightarrow \infty$, further reductions takes place.
As in the case of $\psi(n)$- and $\varphi(n)$-distributions of
random pseudomatrices investigated in [5, Lemma 6.1],
the sum of mixed moments of matricially free random variables 
in the states $\varphi_{j}(n)$ which correspond to a partition 
which is not a non-crossing pair partition is $O(1/\sqrt{n})$.

Moreover, it suffices to take into account the mixed moments of the above type
which are associated with $\pi\in \mathcal{NC}_{m}^{2}(p_1,q_1, \ldots , p_m,q_m)$
for $q=q_m$ since the remaining moments vanish. 
However, these moments are the same for all $(i_1,j_1, \ldots, i_m,j_m)$, for which
$\pi\in \mathcal{NC}_{m}^{2}(i_1,j_1, \ldots , i_m,j_m)$ and $j_m=j$ 
since the variances of $(X_{i,j}(n))$ are block-identical, namely $v_{i,j}(n)=u_{p,q}/n$ 
whenever $(i,j)\in N_p \times N_q$. Thus we can use the matrix elements of $U=(u_{p,q})$ 
to compute our mixed moments and express them as 
$$
\varphi_{j}(n)(X_{i_1,j_1}(n) \ldots X_{i_m,j_m}(n))
=\frac{u_{q}(\pi,f)}{n^{s}},
$$ 
where $f$ is the unique coloring of $\pi\in \mathcal{NC}_{m}^{2}(p_1,q_1, \ldots , p_m,q_m)$
defined by indices $p_{k}$, where $k\in \mathcal{L}(\pi)$, with $q$ coloring the imaginary block. 
Since the coloring $f$ is uniquely determined by $p_1,q_1, \ldots, p_m,q_m$, 
we will use the simplified notation $u_{q}^{*}(\pi)=u_{q}(\pi,f)$.

For such $\pi$ it then remains to enumerate the tuples
$(i_1,j_1,\ldots, i_m,j_m)$, for which
$\pi\in \mathcal{NC}_{m}^{2}(i_1,j_1, \ldots , i_m,j_m)$.
Since the coloring of $\pi$ defined by these tuples is determined by 
independent indices associated with $p_{l(1)}, \ldots, p_{l(s)}$, where
$$
\mathcal{L}(\pi)=\{l(1),l(2), \ldots , l(s)\},
$$
and by $q=q_m$, this enumeration boils down to enumerating
the indices $i_{l(1)}, \ldots , i_{l(s)}$ and $j=j_m$.
An appropriate inductive argument for this fact can be easily 
provided. Instead of giving a formal proof, we refer the reader to Example 2.2 
and to a similar argument given in the context of random matrices 
(see the proof of Theorem 9.1).
Therefore, for given $\pi$, the cardinality of the corresponding set of 
equal mixed moments is of the same order as  
$$
\Theta_{n}(\pi)=n_{q_{m}}n_{p_{l(1)}}\ldots n_{p_{l(s)}}=O(n^{s+1})
$$
as $n\rightarrow \infty$, where $q_m=q$. In fact, the exact cardinality may be slightly smaller 
than $\Theta_{n}(\pi)$ since 
all independent indices $i_{l}$, where $l\in \mathcal{L}(\pi)$, must be different 
in this computation to give a non-zero contribution to the limit.
In the formula for $\Theta_{n}(\pi)$ given above this is not the case 
since the mapping $\mathcal{L}(\pi)\rightarrow [r]$ given by $k\rightarrow p(k)$ may not be 
injective. Nevertheless, we can substitute $\Theta_{n}(\pi)$ for the cardinality
of the considered set when taking the limit $n\rightarrow \infty$ since the 
difference between these two cardinalities is $O(n^{s})$. 

Therefore, from each considered partition we obtain the contribution 
$$
\lim_{n\rightarrow \infty}
\frac{u_{q}^{*}(\pi)\Theta_n(\pi)}{n_{q}n^{s}}=
u_q^{*}(\pi)d_{p_{l(1)}}\ldots d_{p_{l(s)}}
$$
where the division by $n_q$ comes from the normalization of the partial trace
$\psi_{q}(n)$.

Collecting contributions associated with all 
$\pi\in \mathcal{NC}_{m}^{2}(p_1,q_1, \ldots , p_m,q_m)$,
we obtain 
$$
\lim_{n\rightarrow \infty}\psi_{q}(n)(S_{p_1,q_1}(n)\ldots S_{p_m,q_m}(n))=
\sum_{\pi \in \mathcal{NC}_{m}^{2}(p_1,q_1, \ldots , p_m,q_m)}b_{q}^{*}(\pi)
$$
where $b_{q}^{*}(\pi)=b_{q}^{*}(\pi,f)$ corresponds to the unique coloring $f$
of $\pi$ defined by the tuple $(p_1,q_1, \ldots , p_m,q_m)$ and to the matrix $B=DU$.
Finally, the proof that the expression on the right-hand side
is equal to $\psi_{q}(\omega_{p_1,q_1}\ldots \omega_{p_m,q_m})$ is similar to that of Lemma 5.1,
which completes the proof.
\hfill $\blacksquare$

\begin{Corollary}
Under the assumptions of Theorem 7.1, 
the joint $\psi(n)$-distributions of 
blocks and block units converge to the joint $\psi$-distribution of 
the truncated matricially free Gaussian operators 
and truncated units, respectively, as $n\rightarrow \infty$.
\end{Corollary}
{\it Proof.}
If we replace $\psi_{q}(n)$ by $\psi(n)$ at the end of the proof of Theorem 7.1, we need to 
sum over $q\in [r]$ the corresponding right-hand sides 
with weights $d_{q}$, respectively, which result from 
normalizations of partial traces.
Therefore, in order to obtain a combinatorial expression 
for the corresponding mixed moment,
it suffices to replace $b_{q}^{*}(\pi)$ in the above formula by
$b^{*}(\pi)=\sum_{q}d_{q}b_{q}^{*}(\pi)$. 
This gives $\psi(\omega_{p_1,q_1}\ldots \omega_{p_m,q_m})$, which proves our assertion.
\hfill $\blacksquare$

\begin{Corollary}
Under the assumptions of Theorem 7.1, 
the array $(S_{p,q}(n))$ is asymptotically matricially 
free with respect to the array of states $(\psi_{p,q}(n))$
on ${\mathcal A}(n)$ defined by the family of 
normalized partial traces $(\psi_{q}(n))_{1\leq q \leq r}$
as $n\rightarrow \infty$.
\end{Corollary}
{\it Proof.}
It follows from Theorem 7.1 that the $\psi_q(n)$-distribution of 
$(S_{p,q}(n))$ converges 
to the $\psi_q$-distribution of 
$(\omega_{p,q})$ for any $p,q$. 
This, Theorem 7.1 and Proposition 4.2 give our assertion.
\hfill $\blacksquare$

\begin{Theorem}
Under assumptions (A1)-(A4), the joint $\varphi(n)$-distributions of 
blocks and block units converge to the joint 
$\varphi$-distribution of the matricially free Gaussian 
operators and the associated units, respectively, as $n\rightarrow \infty$.
\end{Theorem}
{\it Proof.}
The proof is similar to that of Theorem 7.1 and is based on
the combinatorial arguments of type used in [5, Lemma 6.1] which lead to 
[5, Lemma 6.2].\hfill $\blacksquare$

\begin{Corollary}
Under assumptions (A1)-(A4), the array $(S_{p,q}(n))$ is asymptotically matricially 
free with respect to the array $(\varphi_{p,q}(n))$
of states on ${\mathcal A}(n)$ defined by $\varphi(n)$ and 
the family of normalized partial traces $(\psi_q(n))_{1\leq q \leq r}$ as $n\rightarrow \infty$.
\end{Corollary}
{\it Proof.}
The assertion follows from Theorem 7.2 and Proposition 4.1.
\hfill $\blacksquare$
 
\section{Symmetric matricial freeness}
In order to compare the asymptotics of blocks of random pseudomatrices with 
that of symmetric random blocks, we shall now introduce a symmetric analogue of matricial freeness. 

Roughly speaking, in this concept one replaces ordered pairs by two-element sets
and assumes that the array $({\mathcal A}_{i,j})$ of subalgebras of a unital algebra ${\mathcal A}$ 
contains the diagonal and is symmetric. Moreover, each algebra ${\mathcal A}_{i,j}$ contains an internal unit $1_{i,j}$
which agrees with $1_{j,i}$ for any $(i,j)\in J$. 
By ${\mathcal I}$ we denote the unital algebra
generated by the internal units and we assume that it is commutative.
By $(\varphi_{i,j})$ we denote an array of states on ${\mathcal A}$.  

Instead of sets $\Lambda_{m}$, we shall use their symmetric counterparts, namely subsets of $I^{m}$ of the form
$$
\Pi_{m}=\{(\{i_1,i_2\}, \{i_2,i_3\}, \ldots , \{i_m,i_{m+1}\}): \{i_1,i_2\}\neq \{i_2,i_3\}\neq \ldots \neq \{i_m,i_{m+1}\}\}
$$ 
where $m\in {\mathbb N}$, with their union denoted $\Pi=\bigcup_{m=1}^{\infty}\Pi_m$.
The main difference between `symmetric matricial freeness' and matricial freeness  is that 
in all definitions we have to use $\Pi$ instead of $\Lambda$.

\begin{Definition}
{\rm 
We say that the array $(1_{i,j})$ is a {\it symmetrically matricially free array of units} 
associated with $({\mathcal A}_{i,j})$ and $(\varphi_{i,j})$ if for any diagonal state $\varphi$
it holds that
\begin{enumerate}
\item
$\varphi(u_1au_2)=\varphi(u_1)\varphi(a)\varphi(u_2)$ 
for any $a\in {\mathcal A}$ and $u_1,u_2\in {\mathcal I}$,
\item
if $a_{k}\in {\mathcal A}_{i_k,j_k}\cap {\rm Ker}\varphi_{i_k,j_k}$, where $1<k\leq m$, then
$$
\varphi(a1_{i_1,j_1}a_{2}\ldots a_m)=
\left\{
\begin{array}{cc}
\varphi(aa_{2} \ldots a_m) & {\rm if}\;(\{i_{1},j_{1}\}, \ldots , \{i_{m},j_m\})\in \Pi\\
0 & {\rm otherwise}
\end{array}
\right..
$$
where $a\in {\mathcal A}$ is arbitrary and $\{i_1,j_1\}\neq \ldots \neq \{i_m,j_m\}$.
\end{enumerate}}
\end{Definition}
\begin{Definition}
{\rm We say that a symmetric array  $({\mathcal A}_{i,j})$ is 
{\it symmetrically matricially free} with respect to $({\varphi}_{i,j})$ if
\begin{enumerate}
\item for any $a_{k}\in {\rm Ker}\varphi_{i_k,j_k}\cap {\mathcal A}_{i_k,j_k}$, where $k\in [m]$
and $\{i_1,j_1\}\neq \ldots \neq \{i_m,j_m\}$, and for any diagonal state $\varphi$ it holds that
$$
\varphi(a_1a_2\ldots a_m)=0
$$
\item
$(1_{i,j})$ is a symmetrically matricially free array of units associated with 
$({\mathcal A}_{i,j})$ and $(\varphi_{i,j})$.
\end{enumerate}}
\end{Definition} 
The array of variables $(a_{i,j})$ in a unital algebra ${\mathcal A}$ will be called 
{\it symmetrically matricially free} with respect to $(\varphi_{i,j})$ 
if there exists a symmetrically matricially free array of units 
$(1_{i,j})$ in ${\mathcal A}$ such that the array 
of algebras $({\mathcal A}_{i,j})$, each generated by $a_{i,j}+a_{j,i}$ and $1_{i,j}$, respectively, 
is symmetrically matricially free with respect to $(\varphi_{i,j})$.
The definition of *-matricially free arrays of variables is similar 
to that of *-matricially free arrays.

We shall need the symmetrized operators on the matricially free-boolean Fock space 
${\mathcal F}$ of Section 4, like the $\widehat{\omega}_{i,j}$'s defined as
$$
\widehat{\omega}_{i,j}=
\left\{
\begin{array}{cc}
\omega_{j,j} & {\rm if}\;\;\; i=j\\
\omega_{i,j}+\omega_{j,i}& {\rm if}\;\;\; i\neq j
\end{array}
\right..
$$
In a similar way we define $\widehat{\varsigma}_{i,j}$, 
$\widehat{\zeta}_{i,j}$, $\widehat{\wp}_{i,j}$, etc.
All these arrays are associated with matrix $A=(\alpha_{i,j})$, which is supressed in
the notation. In turn, the symmetrized units on ${\mathcal F}$ and 
${\mathcal F}\ominus {\mathbb C}\Omega$, are of the form 
$$
\widehat{1}_{i,j}=1_{i,j}+1_{j,i}-1_{i,j}1_{j,i}\;\;\; {\rm and}\;\;\;
\widehat{t}_{i,j}=\widehat{1}_{i,j}P,
$$ 
respectively, for any $i,j$, where $(1_{i,j})$ is the array of canonical units
on ${\mathcal F}$ and $P$ stands for the projection onto
${\mathcal F}\ominus {\mathbb C}\Omega$. More explicitly, 
$\widehat{1}_{i,j}$ is the projection onto the subspace of ${\mathcal F}$ spanned 
by tensors which begin with $e_{i,k}$ or $e_{j,k}$ for some $k$, and, in addition,
by $\Omega$ if $i=j$. 

\begin{Proposition}
If the matrix $A=(\alpha_{i,j})$ is symmetric, then 
\begin{enumerate}
\item
the $\psi_j$-distribution of $\,\widehat{\omega}_{i,j}$ is the semicircle law of radius $2\alpha_{i,j}$
for any $(i,j)$,
\item
the array $(\widehat{\omega}_{i,j})$ is symmetrically matricially free with respect to $(\psi_{i,j})$.
\end{enumerate}
\end{Proposition}
{\it Proof.}
If $i=j$, then the first claim is similar to (1) of Proposition 4.1 since the action of 
$\omega_{j,j}$ onto $e_{j,j}$ is similar to the action of $\zeta_{j,j}$ onto the vacuum vector.
If $i\neq j$, then for even and positive $m$, it holds that
$$
\psi_{j}(\widehat{\omega}^{m}_{i,j})
=
\sum_{\epsilon_1\ldots, \epsilon_m\in \{1,*\}}
\sum_{(i_1,j_1), \ldots ,(i_m,j_m)\in \{(i,j),(j,i)\}}
\psi_j\left(\ell^{\,\epsilon_1}(e_{i_1,j_1})
\ldots \ell^{\,\epsilon_{m}}(e_{i_{m},j_{m}})\right)
$$
and there is a bijection between $\mathcal{NC}_{m}^{2}$  and 
products of the form $\ell^{\,\epsilon_1}(e_{i_1,j_1})\ldots \ell^{\,\epsilon_{m}}(e_{i_{m},j_{m}})$
whose action onto $e_{j,j}$ is non-trivial. This bijection is obtained as follows. If $\pi\in \mathcal{NC}_{m}^{2}$
is given and $\{l,r\}$ is a block of $\pi$, where $l<r$, 
then $\epsilon_{l}=*$, $\epsilon_{r}=1$ and $(i_l,j_l)=(i_r,j_r)$ with  
$$
(i_r,j_r)=
\left\{
\begin{array}{cc}
(j_{o(r)},i_{o(r)})&{\rm if}\;\{l,r\}\;{\rm has}\;{\rm an}\;{\rm outer}\;{\rm block}\\
(i,j)& {\rm otherwise}
\end{array}
\right.
$$
where $o(r)$ is the right leg of the nearest outer block of $\{l,r\}$.
It can be seen that this mapping is onto since in order to get a non-trivial action on $e_{j,j}$
of the considered product of operators, the action of each $\ell(e_{i,j})$ 
(corresponding to the right leg of some block) must be followed by the action of
$\ell(e_{j,i})$ (corresponding to the right leg of another block) or $\ell^{*}(e_{i,j})$
(corresponding to the left leg of the same block). In particular, acting first with $\ell(e_{i,j})$
on a vector from ${\mathcal F}$ and then with $\ell(e_{i,j})$ or $\ell^{*}(e_{j,i})$ gives
zero. This completes the proof of (1) since the even moments of the semicircle law of radius $2\alpha_{i,j}$
are given by $M_{m}=\alpha_{i,j}^{m}c_{m/2}$, where $c_{m/2}=|\mathcal{NC}_{m}^{2}|$, $m\in 2{\mathbb N}$,
are Catalan numbers.  
The proof of (2) is similar to that of (2) of Proposition 4.1. 
We show the slightly more general result that the array 
$
\widehat{\mathcal B}_{i,j}={\mathbb C}\langle\widehat{\wp}_{i,j}^{},\widehat{\wp}_{i,j}^{\,*}, 
\widehat{t}_{i,j}^{}\rangle
$
is symmetrically matricially free with respect to $(\psi_{i,j})$, where we use relations
$$
\widehat{\wp}_{i,j}^{\,*}\widehat{\wp}_{i,j}^{}=b_{i,j}^{}\widehat{t}_{i,j}^{},
$$
and normalization $\psi_{j}(\widehat{t}_{i,j})=1$ for any $(i,j)\in J$. 
The details are left to the reader.
\hfill $\blacksquare$\\

For a given partition $[n]=N_{1}\cup N_{2}\cup \ldots \cup N_{r}$ 
of the form considered before, it is useful to introduce the notation
$$
N_{p,q}=(N_{p}\times N_{q})\cup (N_{q}\times N_{p})
$$ 
for sets of pairs which label `symmetric blocks' of random pseudomatrices.
If we let $n\rightarrow \infty$ and assume that blocks grow as before, 
we can obtain the asymptotic behavior of these symmetric blocks.  
\begin{Theorem}
Under the assumptions of Theorem 7.1, the array of symmetric blocks given by
$$
Z_{p,q}(n):=\sum_{(i,j)\in N_{p,q}}X_{i,j}(n)
$$
is asymptotically symmetrically matricially free
with respect to $(\psi_{p,q}(n))$ as $n\rightarrow \infty$.
\end{Theorem}
{\it Proof.}
Observing that
$$
Z_{p,q}(n)=
\left\{
\begin{array}{cl}
S_{p,p}(n) & {\rm if}\;\;p=q\\
S_{p,q}(n)+S_{q,p}(n) & {\rm otherwise}
\end{array}
\right.,
$$ 
it suffices to use Theorem 7.1 and Proposition 8.1 to prove the assertion.
\hfill $\blacksquare$ 

\section{Symmetric random blocks}

We show in this section that the tracial asymptotics of random pseudomatrices is the same as the asymptotics of complex Gaussian random matrices and that this similarity can be carried over to the level of their symmetric random blocks. 

The context for the study of random matrices originated by Voiculescu [10] is 
the following. Let $\mu$ be a probability measure on some measurable space without atoms and let $L=\bigcap_{1\leq p <\infty}L^{p}(\mu)$ be endowed with the state
expectation ${\mathbb E}$ given by integration with respect to $\mu$. 
The *-algebra of $n\times n$ random matrices is 
$M_{n}(L)=L\otimes M_{n}({\mathbb C})$ with the state $\tau(n)={\mathbb E}\otimes 
{\rm tr}(n)$, where ${\rm tr}(n)$ is the normalized trace.

In order to compare the asymptotics of the blocks of 
random matrices with that of random pseudomatrices,
we partition each set $[n]$ into disjoint non-empty intervals as in 
the case of pseudomatrices and we set again $D={\rm diag}(d_1,d_2, \ldots ,d_r)$ 
to be the associated diagonal dimension matrix.

By a {\it complex Gaussian random matrix} we understand a 
matrix $(Y_{i,j}(n))_{1\leq i,j \leq n}$,
in which $Y_{i,j}(n)=\overline{Y_{j,i}(n)}$ for any $i,j,n$ and 
$$
\{{\rm Re}Y_{i,j}(n)|1\leq i \leq j\leq n\}\cup 
\{{\rm Im}Y_{i,j}(n)|1\leq i \leq j\leq n\}
$$ 
is an independent set of Gaussian random variables. 
Its submatrices of the form
$$
T_{p,q}(n)=\sum_{(i,j)\in N_{p,q}}
Y_{i,j}(n)\otimes e_{i,j}(n)
$$
will be called {\it symmetric random blocks}, where $p,q\in [r]$,
and $\{e_{i,j}(n)|1\leq i \leq j \leq n\}$ is a system of matrix units. 

We will study the asymptotics of symmetric random blocks under some natural assumptions.
Our setting is very similar to that in [10,12] except that we will
assume that the variances of $|Y_{i,j}(n)|$, where $i,j\in [n]$ and $n$ is fixed, are block-identical rather than identical.
The case of non-Gaussian random matrices [2] and the corresponding symmetric blocks
can be treated in a similar way.  

In order to find a Hilbert space realization of the limit joint distribution
of symmetric random blocks under $\tau(n)$,
we will use the symmetrized truncated Gaussian operators $\widehat{\omega}_{p,q}$
on ${\mathcal F}$ and the state $\psi=\sum_{j}d_{j}\psi_{j}$ on $B({\mathcal F})$,
where the dependence of $\widehat{\omega}_{p,q}$'s on the matrix $B=DU$ is supressed in the notation.
 
\begin{Theorem}
Let $(Y_{i,j}(n))_{1\leq i,j \leq n}$ be a complex Gaussian random matrix for each $n\in {\mathbb N}$ such that
\begin{enumerate}
\item
${\mathbb E}(Y_{i,j}(n))=0$ for any $i,j,n$,
\item
${\mathbb E}(|Y_{i,j}(n)|^{2})=u_{p,q}/n$ for any $i\in N_p$ and $j\in N_q$ and any $n$,
\end{enumerate}
where $U:=(u_{p,q})\in M_{r}({\mathbb R})$. Then it holds that
$$
\lim_{n\rightarrow \infty}
\tau (n)(T_{p_1,q_1}(n)\ldots T_{p_m,q_m}(n))
=
\psi(\widehat{\omega}_{p_1,q_1}\ldots \widehat{\omega}_{p_m,q_m})\\ 
$$
for any $p_1,q_1, \ldots ,p_m,q_m\in [r]$, where $\widehat{\omega}_{p,q}$'s are
associated with $B:=DU$, where $D$ is the dimension matrix.   
\end{Theorem}
{\it Proof.}
Our proof refers to the original proof of Voiculescu [10,12], 
which is followed by some combinatorial 
arguments referring to non-crossing pair 
partitions (in our approach the variances may vary).
Let $\tau_{q}(n)={\mathbb E}\otimes {\rm tr}_{q}$, 
where ${\rm tr}_{q}(A)=1/n_{q}\sum_{j\in N_{q}}A_{j,j}$
is the normalized partial trace of $A\in M_{n}({\mathbb C})$ 
corresponding to $q\in [r]$. 
We have
$$
\tau_{q}(n)(T_{p_1,q_1}(n)\ldots T_{p_m,q_m}(n))
$$
$$
=
\sum_{(i_1,j_1)\in N_{p_1,q_1}, \ldots , (i_m,j_m)\in N_{p_m,q_m}}
{\mathbb E}(Y_{i_1,j_1}(n)\ldots Y_{i_m,j_m}(n))\;{\rm tr}_{q}(e_{i_1,j_1}\ldots e_{i_m,j_m}).
$$
Since all terms are zero for $m$ odd, throughout the rest of the proof
we assume that $m=2s$ for some $s\in {\mathbb N}$.
  
An individual term of this sum is then non-zero only if 
$$
j_1=i_2, j_2=i_3, \ldots , j_m=i_1\in N_{q}
$$
and there exists a bijection $\gamma: [m]\rightarrow [m]$ such that
$$
\gamma^{2}={\rm id}, \;\;\;\gamma(k)\neq k\;\;{\rm and}\;\;i_{\gamma(k)}=j_k,\; j_{\gamma(k)}=i_k.
$$
Each such term is $O(n^{-s-1})$ as $n\rightarrow \infty$.
For convenience, we restrict our attention to the case when
$$
(i_1,j_1)\in N_{p_1}\times N_{q_1}, \ldots , (i_m,j_m)\in N_{p_m}\times N_{q_m},
$$ 
since the remaining cases can be treated in a similar way, with some $p_{i}$'s 
interchanged with the corresponding $q_{i}$'s.

The number of non-zero terms of the considered type is $\sum_{\gamma}\Theta_{n}(\gamma)$, where 
$\gamma$ runs over the set of permutations of $[m]$ such that 
$$
\gamma^{2}={\rm id},\;\; \gamma(k)\neq k\;\;{\rm and}\;\;p_k=q_{\gamma(k)},\;\;q_{k}=p_{\gamma(k)}
$$
for any $k\in [m]$. If $q=p_1$, we obtain 
$$
\Theta_n(\gamma)={\rm card}\{(i_1,\ldots ,i_m)\in N_{p_1}\times \ldots \times N_{p_m}:
\;i_k=i_{\gamma(k)+1}, \;i_{k+1}=i_{\gamma(k)}\},
$$
where addition is modulo $m$. Denote by $\pi(\gamma)$ the pair-partition 
of $[m]$ defned by $\gamma$. 
 
If $\pi(\gamma)$ is a crossing pair-partition, then $\Theta_{n}(\gamma)\leq O(n^{s})$, as
in the free case, which will make the corresponding mixed moments disappear as $n\rightarrow \infty$ since we still need to multiply these cardinalities by the corresponding expectations, which are $O(1/n^{s})$, and divide them by $n_q$ due to the normalization of the partial trace.

In turn, if $\pi(\gamma)\in \mathcal{NC}_{m}^{2}\setminus \mathcal{NC}_{m}^{2}(\{p_1,q_1\}, \ldots , \{p_m,q_m\})$ or $q\neq p_1$, then again $\Theta_n(\gamma)\leq O(n^{s})$. 
On the other hand, if $\pi(\gamma)\in \mathcal{NC}_{m}^{2}(\{p_1,q_1\}, \ldots , \{p_m,q_m\})$
and $q=p_1$, then the number of independent indices which enter in the computation of $\Theta_n(\gamma)$ is $s+1$. We claim that for such $\gamma$ we have
$$
\Theta_n(\gamma)={\rm card}\{(i_1,i_{r(1)},\ldots , i_{r(s)})\in
N_{p_1}\times N_{p_{r(1)}}\times \ldots \times N_{p_{r(s)}}
\},
$$
where 
$$
\mathcal{R}(\pi(\gamma))=\{r(1),r(2),\ldots ,r(s)\}.
$$ 
In fact, it is not difficult to show that the conditions which define $\Theta_{n}(\gamma)$ can be
reduced to $s+1$ independent indices  
$i_1, i_{r(1)}, \ldots , i_{r(s)}$ which can be interpreted as independent colors assuming arbitrary values
from the corresponding intervals $N_{p_1},N_{p_{r(1)}},\ldots , N_{p_{r(s)}}$, respectively,
with index $i_1$ coloring the imaginary block. 

\begin{figure}
\unitlength=1mm
\special{em:linewidth 0.4pt}
\linethickness{0.4pt}
\begin{picture}(120.00,30.00)(57.00,2.00)
\put(64.00,10.00){\line(0,1){13.50}}
\put(69.00,10.00){\line(0,1){8.00}}
\put(73.00,10.00){$\scriptstyle{{\rm inner}\;{\rm blocks}}$}
\put(93.00,10.00){\line(0,1){8.00}}
\put(97.00,10.00){\line(0,1){8.00}}
\put(121.00,10.00){\line(0,1){8.00}}
\put(102.00,10.00){$\scriptstyle{{\rm inner}\;{\rm blocks}}$}
\put(165.00,10.00){\line(0,1){8.00}}
\put(141.00,10.00){\line(0,1){8.00}}
\put(146.00,10.00){$\scriptstyle{{\rm inner}\;{\rm blocks}}$}
\put(170.00,10.00){\line(0,1){13.50}}
\put(130.00,10.00){$\scriptstyle{\ldots}$}

\put(80.00,20.00){$\scriptstyle{i_{n(2)}}$}
\put(108.00,20.00){$\scriptstyle{i_{n(4)}}$}
\put(115.00,25.50){$\scriptstyle{i_{m(2)}}$}
\put(152.00,20.00){$\scriptstyle{i_{n(2r)}}$}
\put(61.00,7.00){$\scriptscriptstyle{m(1)}$}
\put(67.00,7.00){$\scriptscriptstyle{n(1)}$}
\put(90.00,7.00){$\scriptscriptstyle{n(2)}$}
\put(96.00,7.00){$\scriptscriptstyle{n(3)}$}
\put(119.00,7.00){$\scriptscriptstyle{n(4)}$}
\put(136.00,7.00){$\scriptscriptstyle{n(2r-1)}$}
\put(162.00,7.00){$\scriptscriptstyle{n(2r)}$}
\put(169.00,7.00){$\scriptscriptstyle{m(2)}$}

\put(64.00,23.50){\line(1,0){106.00}}
\put(69.00,18.00){\line(1,0){24.00}}
\put(97.00,18.00){\line(1,0){24.00}}
\put(141.00,18.00){\line(1,0){24.00}}
\end{picture}
\caption{Inductive step involving subpartitions of $\pi(\gamma)$.}
\end{figure}
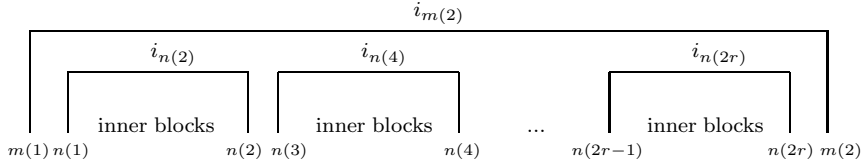

The proof of this fact essentially reduces to the inductive step 
involving subpartitions of $\pi(\gamma)$ of the form shown in Fig.2,
where the covering block can also be interpreted as the imaginary block associated with
$\pi(\gamma)$.
We want to show that the conditions which define 
$\Theta_{n}(\gamma)$ reduce all indices involved here
to those associated with the right legs of the blocks of the subpartition shown above. 
In particular, on the level of blocks of depth one, we have pairings 
$$
\gamma(n(1))=n(2),\gamma(n(3))=n(4), \ldots, \gamma(n(2r-1))=n(2r)
$$ 
which, in view of the condition $i_{k}=i_{\gamma(k)+1}$, lead to
the equation
$$
i_{n(1)}=i_{n(3)}=\ldots =i_{n(2r-1)}=i_{m(2)}
$$
which says that the indices associated with the left legs of blocks of depth one
are equal to the index associated with the right leg of its covering block.
No blocks of depth greater than one lead to new conditions involving 
$i_{m(2)}$ and for that reason the index $i_{m(2)}$ can be used to color the covering block
of the considered subpartition. The same pattern is repeated for blocks of arbitrary depth
which are covered by the same block.
The same holds for blocks of zero depth, where the role of the 
covering block in the above reasoning is played by the imaginary block
(in that case the imaginary block is colored by $i_1$). 

This argument can be viewed as a proof of the inductive step of our claim that 
the conditions involving the indices $i_1, i_2, \ldots , i_m$, used to 
define the numbers $\Theta_{n}(\gamma)$ for the `right' $\gamma$ reduce
to $s+1$ independent indices associated with the right legs of the blocks of $\pi$.
If the covering block is interpreted as the imaginary block, we obtain 
the starting case of the induction. 
This proves our formula for $\Theta_{n}(\gamma)$, from which we easily get 
$$
\Theta_{n}(\gamma)=
n_{p_1}n_{p_{r(1)}}\ldots n_{p_{r(s)}}.
$$
Let us remark here that not all tuples $(i_1,i_{r(1)}, \ldots , i_{r(s)})$
which contribute to $\Theta_{n}(\gamma)$, but only those which are pairwise different,
are actually used in our computation of the asymptotic joint distribution. If any two (or more) 
of these indices coincide, the corresponding contribution from the associated
expectations is $O(n^{-1})$ and therefore becomes irrelevant in the limit.  

For those independent indices which are pairwise different
and correspond to partitions $\pi(\gamma)\in \mathcal{NC}_{m}^{2}(\{p_1,q_1\}, \ldots , \{p_m,q_m\})$,
the corresponding expectation of Gaussian random variables
${\mathbb E}(Y_{i_1,j_1}(n)\ldots Y_{i_m,j_m}(n))$ is the product of
$$
{\mathbb E}(Y_{i_{k},j_{k}}(n)Y_{i_{\gamma(k)},j_{\gamma(k)}}(n))= 
v_{i_{k},j_{k}}(n)
$$
where $\gamma(k)\in \mathcal{R}(\pi(\gamma))$, with $i_{\gamma(k)}=j_k$ and 
$j_{\gamma(k)}=i_{k}$. 
Now, since $i_{k}=i_{o(k)}$, 
the above expectation is equal to $v_{i,j}(n)$,
where $i=i_{o(k)}$ and $j=i_{\gamma(k)}$.
In other words, $j$ is the color assigned to block $\{k,\gamma(k)\}$ 
since we have shown before that we color the blocks 
with indices which correspond to their right legs
and $i$ is the color assigned to its covering block.
This includes the case
of blocks which do not have an outer block, then 
the role of the latter is played by the conditional 
block colored  by $i_1=j_m$.

Now, using our assumption on the variance matrices,
$v_{i,j}(n)=v_{j,i}(n)=u_{q,p}/n$ whenever $i\in N_p, j\in N_q$.
Therefore, taking into account all pairings and using the definition of 
symbol $u_j(\pi,f)$ of Section 2, we obtain
$$
{\mathbb E}(Y_{i_1,j_1}(n)\ldots Y_{i_m,j_m}(n))=
\frac{u_{p_1}(\pi, f)}{n^{s}}
$$ 
for the relevant Gaussian expectations of variables 
associated with $\pi(\gamma)$,
where $f\in F_{r}(\pi(\gamma))$ 
is the coloring of the associated $\pi(\gamma)$ defined by 
indices $p_{r(1)}, p_{r(2)}, \ldots , p_{r(s)}$ assigned to the blocks
of $\pi(\gamma)$, with color $p_{1}$ assigned to the imaginary block.

Collecting expectations corresponding to all $\pi\in \mathcal{NC}_{m}^{2}$
and taking into account that 
$$
\lim_{n\rightarrow \infty}\frac{\Theta_{n}(\gamma)}{n^{s}}=
d_{p_{1}}d_{p_{r(1)}}\ldots d_{p_{r(s)}}, 
$$
we obtain
$$
\lim_{n\rightarrow \infty}\tau_{q}(n)(T_{p_1,q_1}(n)\ldots T_{p_m,q_m}(n))=
\sum_{\pi \in 
\mathcal{NC}_{m}^{2}(\{p_1,q_1\},\ldots \{p_m,q_m\})}
\widehat{b}_{q}(\pi)
$$
where $\widehat{b}_q(\pi)=\sum_{f}b_{q}(\pi,f)$ is the sum over all admissible colorings
of $\pi$ and each $b_{q}(\pi,f)$ coresponds to the matrix $B=DU$ and
to the color $q$ of the imaginary block.

If we replace $\tau_{q}(n)$ by 
$\tau(n)$ in the above expression, 
we just need to sum over $q\in [r]$ the corresponding right-hand sides 
with weights $d_{q}$, respectively, which result from 
normalizations of partial traces.
Therefore, in order to obtain the required expression for the corresponding mixed moment,
we need to replace $\widehat{b}_{q}(\pi)$ by $\widehat{b}(\pi)=\sum_{q}d_{q}\widehat{b}_{q}(\pi)$, 
which completes the proof of the combinatorial formula.

It remains to show that the mixed moment of the 
symmetrized Gaussian operators in the state $\psi$ 
gives the apropriate sum of combinatorial expressions obtained above.
This, in turn, is similar to the proof of Lemma 5.1. Therefore, our proof is completed.
\hfill $\blacksquare$ \\

\begin{Remark}
{\rm In this context, let us remark that the array of symmetric random blocks
$(T_{p,q}(n))$ is not symmetrically matricially free with respect
to $(\tau_{p,q}(n))$ for finite $n$. In order to see this, it suffices to compute
some simple examples of mixed moments of the $T_{p,q}(n)$'s
associated with crossing partitions.
For instance, if $n=3$ and the blocks are one-dimensional, then 
$\tau_{1}(T_{1,2}T_{2,3}T_{3,1}T_{1,2}T_{2,3}T_{3,1})\neq 0$,
where $\tau_1=\tau_1(3)$ and $T_{p,q}=T_{p,q}(3)$ for any $1\leq p,q \leq 3$,
which shows that condition (1) of Definition 8.2 is not satisfied.}
\end{Remark}

However, in view of Proposition 8.1 and Theorem 9.1, we can expect that 
symmetric random blocks are symmetrically matricially free 
as $n\rightarrow \infty$. The notion of asymptotic symmetric matricial freeness
is analogous to that of asymptotic matricial freeness (see Definition 7.1).
It remains to define appropriate arrays of 
block units. Thus, let
$$
\widehat{\bf 1}_{p,q}(n)=\sum_{j\in N_{p}\cup N_{q}}1\otimes e_{j,j}(n)
$$
for any $p,q\in [r]$ and each $n\in {\mathbb N}$. It is easy to see that
$\widehat{\bf 1}_{p,q}(n)$ is an internal unit in the algebra generated by
$T_{p,q}(n)$ and $\widehat{\bf 1}_{p,q}(n)$ for any given $p,q\in [r]$ and $n\in {\mathbb N}$.
\begin{Theorem}
Under the assumptions of Theorem 9.1, let $(\tau_{p,q}(n))$ be
the $r$-dim\-en\-sio\-nal array of states on the algebra 
${\mathcal A}(n)=L\otimes M_{n}({\mathbb C})$ defined by
the normalized partial traces $\tau_{q}(n)$, where $q\in [r]$, for any
$n\in {\mathbb N}$. Then the array of symmetric random blocks $(T_{p,q}(n))$ is asymptotically
symmetrically matricially free with respect to $(\tau_{p,q}(n))$.
\end{Theorem}
{\it Proof.}
In view of Proposition 8.1, it suffices to show that 
the $\tau_m(n)$-distributions of blocks $T_{p,q}(n)$ and units 
$\widehat{\bf 1}_{p,q}(n)$ converge to the $\psi_m$-distribution of 
the symmetrized Gaussians $\widehat{\omega}_{p,q}$ and summetrized units 
$\widehat{1}_{p,q}$ as $n\rightarrow \infty$ for any given $m$. It suffices to verify that 
$\widehat{\bf 1}_{p,q}(n)$ leaves invariant any vector of the form 
$$
w_1(T_{p_1,q_1}(n))\ldots w_k(T_{p_k,q_k}(n))e_m 
$$
where the product of polynomials $w_1, \ldots ,w_k$ 
is in the symmetrically matricially free kernel form with respect to 
the array $(\tau_{p,q}(n))$ and $\{p,q\}\cap \{p_1,q_1\}\neq \emptyset$
and kills all other vectors. This is quite easy 
to observe since any vector of this type is a linear combination 
of the form $\sum_{i\in N_p}\alpha_i a_i\otimes e_i+ \sum_{j\in N_q}\beta_{j} 
b_j\otimes e_j$, where $a_i,b_j\in L$ and $\alpha_i,\beta_j\in {\mathbb C}$ for any $i,j$,
and thus it is left invariant by $\widehat{\bf 1}_{p,q}(n)$,
where $(e_j)$ is the canonical basis in ${\mathbb C}^{n}$.
\hfill $\blacksquare$

\section{Asymptotic freeness and asymptotic monotone independence}

If we make additional assumptions on the variance matrices $(v_{i,j}(n))$ 
of matricially free arrays $(X_{i,j}(n))$, we obtain
asymptotic freeness which refers to the rows of 
random pseudomatrices.
The proposition given below can also be viewed as an operatorial
version of [5, Proposition 8.2]. 

\begin{Proposition}
If the arrays $(X_{i,j}(n))$ of Theorem 3.2 are square and have identical variances within blocks and rows, then the sums
$$
S_{p}(n):=\sum_{q=1}^{r}S_{p,q}(n), \;\;\; where \;\;\;1\leq p \leq r,
$$ 
are asymptotically free with respect to both $\varphi(n)$ and $\psi(n)$ as $n\rightarrow \infty$.
\end{Proposition}
{\it Proof.}
In view of Theorems 7.1-7.2, it suffices to show that the variables $\zeta_1, \ldots , \zeta_r$
are free with respect to $\varphi$, and that $\omega_1, \ldots , \omega_r$
are free with respect to $\psi$, where
$$
\zeta_{p}=\sum_{q=1}^{r}\zeta_{p,q}\;\;\;{\rm and}\;\;\;  \omega_{p}=\sum_{q=1}^{r}\omega_{p,q},
$$
with the notations of Section 4. 
If $\{e_1,e_{2}, \ldots , e_r\}$ 
is an orthonormal set of vectors, we have the natural isomorphism
$$
\kappa:{\mathcal F}\rightarrow {\mathcal F}(\bigoplus_{q}{\mathbb C}e_{q})
$$
of Remark 4.1 and then, since $\alpha_{p,q}=\alpha_{p,p}$ for any $p,q \in [r]$, 
we have
$$
\zeta_{p}=\alpha_{p,p}\kappa^{*}\omega(e_{p})\kappa
$$
where $\omega(e_{p})=\ell(e_{p})+\ell^{*}(e_{p})$ is the canonical
free Gaussian operator on ${\mathcal F}(\bigoplus_{q}{\mathbb C}e_{q})$ for any $p$.
Let us remark that $\omega_p$ and $\omega(e_p)$ denote different operators in our notation.
Since $\omega(e_{1}), \ldots , \omega(e_{r})$ 
are free with respect to the vacuum state on $B({\mathcal F}(\bigoplus_{q}{\mathbb C}e_{q}))$,
the operators $\zeta_{1}, \ldots , \zeta_r$ are free with respect
to the vacuum state on $B({\mathcal F})$.
Similarly,
$$
\omega_{p}=\alpha_{p,p}\gamma^{*}\omega(e_{p})\gamma
$$
for any $p\in [r]$, where 
$\gamma=P\circ \kappa$ and $P$ is the canonical projection from 
${\mathcal F}(\bigoplus_{q}{\mathbb C}e_{q})$ onto the orthocomplement of ${\mathbb C}\Omega$.
Now, since $P\omega(e_{1})P, \ldots , P\omega(e_r)P$ are 
free with respect to $\psi_{q}(.)=\langle .e_{q}, e_{q}\rangle$ for any $q$,
the operators $\omega_{1},\ldots , \omega_{r}$ are free with respect to
$\psi_{q}$ for any $q$. Moreover, the $\psi_p$-distribution of $\omega_{p}$ does not 
depend on $p$ since the variances are assumed to be identical in each row.
Therefore, $\omega_{1},\ldots , \omega_{r}$ are free with respect to $\psi$. 
This completes the proof.
\hfill $\blacksquare$\\

Analogous results hold for block-triangular arrays and lead to monotone
independence (lower-block-triangular arrays) 
and anti-monotone independence (upper-block-triangular arrays).
One has to remember that order is important for these notions of independence
and therefore in that case we use sequences of variables instead of families. 
We formulate the result only for 
lower-triangular arrays since the case of 
upper-triangular arrays is completely analogous.
This result can be viewed as an operatorial version of [5, Proposition 8.3].

\begin{Proposition}
If the arrays $(X_{i,j}(n))$ of Theorem 3.2 are lower-block-triangular 
and have identical variances within blocks and rows, then the sums
$$
S_{p}(n):=\sum_{q=1}^{p}S_{p,q}(n), \;\;\; where \;\;\;1 \leq p \leq r,
$$ 
are asymptotically monotone independent with respect 
to $\varphi(n)$ as $n \rightarrow \infty$.
\end{Proposition}
{\it Proof.}
Let $\zeta_{p}=\sum_{q=1}^{p}\zeta_{p,q}$,
where $p \in [r]$. In view of Theorem 7.1, it suffices to show that
the operators $\zeta_{1}, \ldots , \zeta_{r}$ are monotone independent with respect to $\varphi$. 
Let ${\mathcal F}_{1}$ be the subspace of ${\mathcal F}(\bigoplus_{q}{\mathbb C}e_{q})$ of the form 
$$
\mathcal{F}_{1}={\mathbb C}\xi\oplus\bigoplus_{m=1}^{r}\bigoplus_{p_1>\ldots >p_m}\bigoplus_{n_1, \ldots , n_m\in {\mathbb N}}
{\mathcal H}_{p_1}^{\otimes n_1}\otimes \ldots \otimes {\mathcal H}_{p_m}^{\otimes n_m}
$$
where ${\mathcal H}_{q}={\mathbb C}e_{q}$ for any $q\in [r]$ and 
$\{e_1, e_{2},\ldots , e_r\}$ is an orthonormal basis in some Hilbert space and 
denote by $Q:{\mathcal F}(\bigoplus_{q}{\mathbb C}e_{q})\rightarrow \mathcal{F}_{1}$
the corresponding canonical projection. Observe that we have
$$
\zeta_{p}=\alpha_{p,p}\beta^{*}\omega(e_{p})\beta
$$
where $\beta= Q \circ \kappa$ and $\kappa$ is the same as in the proof of 
Proposition 10.1. 
Since the operators $Q\omega(e_1)Q, \ldots, Q\omega(e_r)Q$
are monotone independent with respect to the vacuum state on the free Fock space, 
the operators $\zeta_{1}, \ldots , \zeta_r$ are monotone independent with respect to $\varphi$.
This completes the proof.
\hfill $\blacksquare$\\

Finally, we obtain asymptotic boolean independence of 
blocks of block-diagonal pseudomatrices.
 
\begin{Proposition}
If the arrays $(X_{i,j}(n))$ of Theorem 3.2 are block-diagonal 
and have identical variances within blocks, then the sums
$S_{p,p}(n)$, where $1 \leq p \leq r$, 
are asymptotically boolean independent with respect 
to $\varphi(n)$ as $n \rightarrow \infty$.
\end{Proposition}
{\it Proof.}
This is a simple consequence of Theorem 7.2 and the fact that the family
$(\zeta_{p,p})_{1\leq p \leq r}$ is boolean independent with respect to the vacuum state $\varphi$ on ${\mathcal F}$.
\hfill $\blacksquare$


\begin{thebibliography}{99}
\bibitem{[1]}
Ph. Biane, Processes with free increments, {\it Math. Z.} {\bf 227} (1998), 143-174.
\bibitem{[2]}
K. Dykema, On certain free product factors via an extended matrix model, 
{\it J. Funct. Anal.} {\bf 112} (1993), 31-60.
\bibitem{[3]}
R. Lenczewski, Decompositions of the free additive convolution, 
{\it J. Funct. Anal.} {\bf 246} (2007), 330-365.
\bibitem{[4]}
R. Lenczewski, Operators related to subordination for free multiplicative convolutions,
{\it Indiana Univ. Math. J.}, {\bf 57} (2008), 1055-1103.
\bibitem{[5]}
R. Lenczewski, Matricially free random variables,  arXiv:0812.0488v1 [math.OA], 2008.
\bibitem{[6]}
N. Muraki, Monotonic independence, monotonic central limit theorem and monotonic law of small numbers,
{\it Infin. Dimens. Anal. Quantum Probab. Relat. Top.} {\bf 4} (2001), 39-58.
\bibitem{[7]}
D. Shlyakhtenko, Random Gaussian band matrices and freeness with amalgamation, {\it Int. Math. Res. Notices}
{\bf 20} (1996), 1013-1025.
\bibitem{[8]}
D. Voiculescu, Symmetries of some reduced free product $C^{*}$-algebras, Operator Algebras and Their Connections with Topology and Ergodic Theory, Lecture Notes in Mathematics, Vol. 1132, Springer Verlag, 1985, pp. 556-588.
\bibitem{[9]}
D. Voiculescu, Lectures on free probability theory, {\it Lectures on probability theory and statistics}
(Saint-Flour, 1998), 279-349, Lecture Notes in Math. 1738, Springer, Berlin, 2000.
\bibitem{[10]}
D. Voiculescu, Limit laws for random matrices and free products,
{\it Invent. Math.} {\bf 104}(1991), 201-220.
\bibitem{[11]}
D. Voiculescu, The analogues of entropy and of Fisher's information measure in free probability theory, I,
{\it Commun. Math. Phys.} {\bf 155} (1993), 71-92.
\bibitem{[12]}
D. Voiculescu , K. Dykema, A. Nica, {\it Free random variables}, CRM Monograph
Series, No.1, A.M.S., Providence, 1992.
\bibitem{[13]}
E. Wigner, On the distribution of the roots of certain symmetric matrices, 
{\it Ann. Math.} {\bf 67} (1958), 325-327.
\end{thebibliography}
\end{document}